\documentclass[11pt]{article}
\usepackage{amsfonts}

\usepackage{amsmath}
\usepackage{latexsym,amsfonts,amssymb}

\newcommand{\epsln}{\varepsilon}

\newcommand{\bfe}{{\mathbf e}}

\newcommand{\bfv}{{\mathbf v}}
\newcommand{\bfw}{{\mathbf w}}
\newcommand{\bfg}{{\mathbf g}}
\newcommand{\bfh}{{\mathbf h}}

\newcommand{\bfn}{{\mathbf n}}

\newcommand{\bfI}{{\mathbf I}}
\newcommand{\bfL}{{\mathbf L}}

\newcommand{\bfT}{{\mathbf T}}
\newcommand{\bfV}{{\mathbf V}}
\newcommand{\bfW}{{\mathbf W}}

\newcommand{\hs}{\hspace{0.3cm}}

\newcommand{\vsss}{\vspace{.05in}}

\begin{document}

\title{Asymptotic Stability of Stationary Solutions of a Free
  Boundary Problem Modeling the Growth of Tumors with Fluid Tissues}
\author{Junde Wu$^{\dag}$\ \ \ and\ \ \ Shangbin Cui$^{\ddag}$\\[0.2cm]
  {\small $^{\dag}$ Department of Mathematics, Sun Yat-Sen University,
  Guangzhou, Guangdong 510275,}\\ [-0.1cm]
  {\small People's Republic of China. E-mail:\,wjdmath@yahoo.com.cn}
  \\[0.1cm]
  {\small $^{\ddag}$ Institute of Mathematics, Sun Yat-Sen University,
  Guangzhou, Guangdong 510275,}\\ [-0.1cm]
  {\small People's Republic of China. E-mail:\,cuisb3@yahoo.com.cn}}
\date{}
 \maketitle

\begin{abstract}
  This paper aims at proving asymptotic stability of the radial stationary
  solution of a free boundary problem modeling the growth of nonnecrotic
  tumors with fluid-like tissues. In a previous paper we considered the
  case where the nutrient concentration $\sigma$ satisfies the stationary
  diffusion equation $\Delta\sigma=f(\sigma)$, and proved that there exists
  a threshold value $\gamma_*>0$ for the surface tension coefficient $\gamma$,
  such that the radial stationary solution is asymptotically stable in case
  $\gamma>\gamma_*$, while unstable in case $\gamma<\gamma_*$. In this paper
  we extend this result to the case where $\sigma$ satisfies the non-stationary
  diffusion equation $\epsln\partial_t\sigma=\Delta\sigma-f(\sigma)$. We prove
  that for the same threshold value $\gamma_*$ as above, for every $\gamma>
  \gamma_*$ there is a corresponding constant $\epsln_0(\gamma)>0$ such that
  for any $0<\epsln<\epsln_0(\gamma)$ the radial stationary solution is
  asymptotically stable with respect to small enough non-radial perturbations,
  while for $0<\gamma<\gamma_*$ and $\epsln$ sufficiently small it is unstable
  under non-radial perturbations.
\medskip

\medskip


\end{abstract}

\section{Introduction}
\setcounter{equation}{0}

  This paper is concerned with the following free boundary problem modelling
  the growth of tumors with fluid-like tissues:
\begin{equation}
   \epsln\partial_t\sigma=\Delta\sigma-f(\sigma) \quad\mbox{in}\;\;
    \Omega(t),\;\;t>0,
\end{equation}
\begin{equation}
  \nabla\cdot\bfv=g(\sigma) \quad\mbox{in}\;\;
  \Omega(t),\;\;t\ge0,
\end{equation}
\begin{equation}
   -\nu\Delta\bfv+\nabla p-{\nu\over3}\nabla
   (\nabla\cdot\bfv)=0 \quad\mbox{in} \;\;\Omega(t),\;\;t\ge0,
\end{equation}
\begin{equation}
   \sigma=\bar\sigma \quad\mbox{on}\;\;\partial\Omega(t),\;\;t\ge0,
\end{equation}
\begin{equation}
  \bfT({\bfv},p){\bfn}=-\gamma\kappa{\bfn}
  \quad\mbox{on}\;\;\partial\Omega(t),\;\;t\ge0,
\end{equation}
\begin{equation}
   V_n=\bfv\cdot\bfn\quad\mbox{on}\;\;\partial\Omega(t),\;\;t\ge0,
\end{equation}
\begin{equation}
   \int_{\Omega(t)}\bfv\;dx=0, \quad t>0,
\end{equation}
\begin{equation}
   \int_{\Omega(t)}\bfv\times{x}d\;x=0,\quad t>0,
\end{equation}
\begin{equation}
   \sigma(0,x)=\sigma_0(x)\quad\mbox{for}\;\; x\in\Omega_0,
\end{equation}
\begin{equation}
   \Omega(0)=\Omega_0,
\end{equation}
  where $\sigma=\sigma(t,x)$, $\bfv=\bfv(t,x)$ ($=(v_1(t,x),v_2(t,x),v_3(t,x))$)
  and $p=p(t,x)$ are unknown functions representing the concentration of nutrient,
  the velocity of fluid and the internal pressure, respectively,
  $f$ and $g$ are given functions representing the nutrient consumption rate and
  tumor cell proliferation rate, respectively, and
  $\Omega(t)$ is an a priori unknown bounded domain in $\Bbb R^3$ representing the
  region occupied by the tumor at time $t$. Besides, $\epsln$, $\nu$, $\bar\sigma$ and
  $\gamma$ are positive constants, among which $\epsln$ is the ratio between typical,
  $\nu$ is the viscosity coefficient of the tumor tissue, $\gamma$
  is the surface tension coefficient of the tumor surface, and $\bar\sigma$
  is the concentration of nutrient in tumor's host tissues, $\kappa$, $V_n$
  and $\bfn$ denote the mean curvature, the normal velocity and the unit outward
  normal, respectively, of the tumor surface $\partial\Omega(t)$, and
  $\bfT({\bfv},p)$ denotes the stress tensor, i.e.,
\begin{equation}
  \bfT({\bfv},p)=\nu\big[\nabla\otimes\bfv+(\nabla\otimes\bfv)^T\big]-
  (p+{2\nu\over3}\nabla\cdot\bfv)\bfI,
\end{equation}
  where $\bfI$ denotes the unit tensor. We note that the sign of the mean
  curvature $\kappa$ is defined such that it is nonnegative for convex
  hyper-surfaces. Without loss of generality, later on we assume that $\nu=1$
  and $\bar{\sigma}=1$. Note that the general situation can be easily reduced
  into this special situation by rescaling.
  As in \cite{WuCui}, throughout this paper we assume that $f$ and $g$  are
  generic smooth functions satisfying the following assumptions:

\medskip
  (A1) $f\in C^\infty[0,\infty)$, $f'(\sigma)>0$ for $\sigma\ge0$
  and $f(0)=0$.

  (A2) $g\in C^\infty[0,\infty)$, $g'(\sigma)>0$ for $\sigma\ge0$
  and $g(\tilde\sigma)=0$ for some $\tilde\sigma>0$,

  (A3) $\tilde\sigma<\bar\sigma$.
\medskip

  The above problem is a simplified form of the tumor models proposed by
  Franks et al in literatures \cite{Franks1}--\cite{Franks4}, which mimic the
  early stages of the growth of ductal carcinoma in breast, and was first
  studied by Friedman in \cite{Fried1}. Local well-posedness of this problem
  in H\"{o}lder spaces has been established by Friedman \cite{Fried1} in a
  more general setting. Moreover, in \cite{Fried1} it is also proved that,
  for the special case $f(\sigma)=\lambda\sigma$ and $g(\sigma)=\mu(\sigma-
  \tilde{\sigma})$, the problem (1.1)--(1.10) has a unique radially symmetric
  stationary solution $(\sigma_s,\bfv_s,p_s,\Omega_s)$. In \cite{FriedHu2}
  Friedman and Hu proved that there exists a threshold value $(\mu/\gamma)_*$
  such that in the case $\mu/\gamma<(\mu/\gamma)_*$ this radial stationary
  solution is {\em linearly asymptotically stable}, i.e. the trivial solution
  of the linearization at $(\sigma_s,\bfv_s,p_s,\Omega_s)$ of the original
  problem is asymptotically stable, and in the case $\mu/\gamma>(\mu/\gamma)_*$
  this stationary solution is unstable. However, whether or not in the case
  $\mu/\gamma<(\mu/\gamma)_*$ this stationary solution is {\em asymptotically
  stable}, namely, whether or not $(\sigma_s,\bfv_s,p_s,\Omega_s)$ is
  asymptotically stable under arbitrary sufficiently small non-radial
  perturbations, which is the Problem 3 of \cite{Fried1} (see also the Open
  Problem $(i)$ in Section 2 of \cite{Fried2}), was not answered by these
  mentioned literatures.

  In a previous work (see \cite{WuCui}) we studied the above problem for the
  model simplified from (1.1)--(1.10) by taking $\epsln=0$, and proved that
  there exists a threshold value $\gamma_*>0$ for the surface tension coefficient
  $\gamma$, such that in the case $\gamma>\gamma_*$ the radial stationary
  solution is asymptotically stable with respect to small enough non-radial
  perturbations, while in case $\gamma<\gamma_*$ this stationary solution is
  unstable under non-radial perturbations. The aim of the present work is to
  extend this result to the case where $\epsln$ is non-vanishing but small. We
  shall prove that for the same threshold value $\gamma_*$ as above, for every
  $\gamma>\gamma_*$ there is a corresponding constant $\epsln_0(\gamma)>0$ such
  that for any $0<\epsln<\epsln_0(\gamma)$ the radial stationary solution is
  asymptotically stable with respect to small enough non-radial perturbations,
  while for $0<\gamma<\gamma_*$ and $\epsln$ sufficiently small it is unstable
  under non-radial perturbations. To give a precise statement of our main
  result, let us first introduce some notations.

  As in \cite{WuCui}, we denote by $(\sigma_s,\bfv_s,p_s,\Omega_s)$ the
  unique radial stationary solution of (1.1)--(1.8), i.e., $\Omega_s=\{x\in
  \Bbb R^3:\,|x|<R_s\}$ and
$$
  \sigma_s=\sigma_s(r), \quad \bfv_s=v_s(r)\frac{x}{r}, \quad
  p_s=p_s(r) \quad \mbox{for}\;\; x\in\Omega_s,
$$
  and for any $x_0\in \Bbb R^3$ we denote by
  $(\sigma_{[x_0]},\bfv_{[x_0]},p_{[x_0]},\Omega_{[x_0]})$ the stationary
  solution of (1.1)--(1.8) obtained by the coordinate translation $x\to x+x_0$
  of the stationary solution $(\sigma_s,\bfv_s,p_s,\Omega_s)$. Given $\rho\in
  C^1(\partial\Omega_s)$ with $\|\rho\|_{C^1(\partial\Omega_s)}$ sufficiently
  small, we denote by $\Omega_\rho$ the domain enclosed by the hypersurface
  $r=R_s+\rho(\xi)$, where $\xi\in\partial\Omega_s$. Since we shall only be
  concerned with small perturbations of the stationary solution $(\sigma_s,
  \bfv_s,p_s,\Omega_s)$, there exist functions $\rho(t)$ ($=\rho(\xi,t)$)
  and $\rho_0$ ($=\rho_0(\xi)$) on $\partial\Omega_s$ such that $\Omega(t)=
  \Omega_{\rho(t)}$ and $\Omega_0=\Omega_{\rho_0}$. Using these notations, the
  initial condition (1.10) can be rewritten as follows:
\begin{equation}
  \rho(\xi,0)=\rho_0(\xi) \quad \mbox{for}\;\; \xi\in\partial\Omega_s.
\end{equation}
  The solution $(\sigma,{\bfv},p,\Omega)$ of the problem (1.1)--(1.9) will be
  correspondingly rewritten as $(\sigma,{\bfv},p,\rho)$, and the radially
  symmetric stationary solution $(\sigma_s,\bfv_s,p_s,\Omega_s)$ will be
  re-denoted as $(\sigma_s,\bfv_s,p_s,0)$.

  The  main result of this paper is the following theorem:
\medskip

  {\bf Theorem 1.1}\ \ {\em Assume that Assumptions $(A1)$--$(A3)$ hold. For
  given $m\in {\Bbb N}$, $m\geq 3$, and $0<\theta<1$, we have the following
  assertion: There exists a positive threshold value $\gamma_*$ such that for
  any $\gamma>\gamma_*$, the radially symmetric stationary solution $(\sigma_s,
  \bfv_s,p_s,0)$ is asymptotically stable for small $\epsln$ in the following
  sense: There exists $\epsln_0>0$ such that for any $0<\epsln<\epsln_0$ there
  exists a corresponding constant $\epsilon>0$, such that for any $\rho_0\in
  C^{m+\theta}(\partial\Omega_s)$ satisfying $||\rho_0||_{C^{m-1+\theta}
  (\partial\Omega_s)}<\epsilon$, the problem $(1.1)$--$(1.9)$ has a unique
  solution $(\sigma,\bfv,p,\rho)$ for all $t\geq 0$, and there are positive
  constants $\omega$, $K$ independent of the initial data and a point
  $x_0\in\Bbb R^3$ uniquely determined by the initial data, such that
  the following holds for all $t\geq 1$:
\begin{eqnarray}
  &&||\sigma(\cdot,t)-\sigma_{[x_0]}||_{C^{m+\theta}(\Omega(t))}+
  ||\bfv(\cdot,t)-\bfv_{[x_0]}||_{C^{m-1+\theta}(\Omega(t))} \nonumber
  \\
  &+&||p(\cdot,t)-p_{[x_0]}||_{C^{m-2+\theta}(\Omega(t))}
  +||\rho(\cdot,t)-\rho_{[x_0]}||_{C^{m+\theta}(\partial\Omega_s)}
  \le K e^{-\omega t}.
\end{eqnarray}
  For $\gamma<\gamma_*$ and $\epsln$ sufficiently small, the stationary
  solution $(\sigma_s,{\bfv}_s,p_s,0)$ is unstable. $\quad$$\Box$}
\medskip

  As in \cite{WuCui} we shall use a functional approach to prove this result,
  namely, we shall first reduce the problem (1.1)--(1.10) into a differential
  equation in a Banach space, and next use the geometric theory for differential
  equations in Banach spaces to study the  asymptotic behavior of the reduced
  equation. However, unlike the case $\epsln=0$ in which the reduced equation
  is a scalar (first-order) nonlinear parabolic pseudo-differential equation on
  the compact manifold $\Bbb S^2$ which does not have a boundary, in the present
  case $\epsln\neq 0$ the reduced equation is a system of equations, one of
  which has a similar feature as the equation in the case $\epsln=0$, while the
  other of which is defined on the domain $\Omega_s$ complemented with a
  Dirichlet boundary condition.
  This determines that in the present case we are forced to deal with a number
  of new difficulties. The first difficulty lies in computation of the spectrum
  of the linearized operator, because we now encounter a matrix operator which
  is not of the diagonal form. To overcome this difficulty we shall use a
  technique developed in \cite{Cui2} to show that the linearized operator is
  similar to a small perturbation of a matric operator possessing a triangular
  structure. The second difficulty is caused by the Dirichlet boundary condition,
  which determines that we cannot find a suitable continuous interpolation space
  as our working space to make the center manifold analysis. More precisely, as
  in the case $\epsln=0$, $0$ is an eigenvalue of the linearized operator, so
  that the standard linearized stability principle does not apply. In the case
  $\epsln=0$, this difficulty is overcome with the aid of the center manifold
  analysis technique developed in \cite{EscSim}. Since this technique requires
  that the working space must be a continuous interpolation space, it fails to
  apply to the present case $\epsln\neq 0$. To overcome this difficulty we
  shall use the idea of Lie group action developed in \cite{Cui2} and apply
  Theorem 2.1 of \cite{Cui2} to solve this problem.

  The structure of the rest part is as follows. In Section 2 we first
  convert the problem into an equivalent initial-boundary value problem
  on a fixed domain by using {\em Hanzawa transformation}, and next we further
  reduce it into a differential equation in a Banach space. In Section 3 we
  study the linearization of (1.1)--(1.8) at the radial stationary solution, and
  compute the spectrum of the linearized operator. In the last section we give
  the proof of Theorem 1.2.

\section{Reduction of the problem}
\setcounter{equation}{0}

  In this section we reduce the problem (1.1)--(1.10) into a differential
  equation in a Banach space. For simplicity of the notation, later on
  we always assume that $R_s=1$. Note that this assumption is reasonable
  because the case $R_s\neq 1$ can be easily reduced into this case after
  a rescaling. It follows that
$$
  \Omega_s=\Bbb B^3=\{x\in\Bbb R^3:|x|<1\}\quad \mbox{and} \quad
  \partial\Omega_s=\Bbb S^2=\{x\in\Bbb R^3:|x|=1\}.
$$

  Let $m$ and $\theta$ be as in Theorem 1.1. We introduce an operator
  $\Pi\in L(C^{m+\theta}(\Bbb S^2),C^{m+\theta}(\overline{\Bbb B}^3))$
  as follows: Given $\rho\in C^{m+\theta}(\Bbb S^2)$, we define
  $u=\Pi(\rho)\in C^{m+\theta}(\overline{\Bbb B}^3)$ to be the
  unique solution of the boundary value problem
\begin{equation}
  \Delta u=0\quad \mbox{in}\;\;\Bbb B^3,\qquad u=\rho
  \quad\mbox{on}\;\;\Bbb S^2.
\end{equation}
  It is clear that $\Pi\in L(C^{m+\theta}(\Bbb S^2),C^{m+\theta}(\overline{\Bbb B}^3)$,
  and $\Pi$ is a right inverse of the trace operator, i.e., we
  have $\mbox{tr}_{\Bbb S^2}\big(\Pi(u))\big)=u$ for any $u\in C^{m+\theta}
  (\Bbb S^2)$. Let $E\in L(C^{m+\theta}(\overline{\Bbb B}^3),BUC^{m+\theta}
  ({\Bbb R}^3))$ be an extension operator, i.e., $E$ has the property that
  $E(u)(x)=u(x)$ for any $u\in C^{m+\theta}(\overline{\Bbb B}^3)$ and $x\in
  \overline{\Bbb B}^3$. Here $BUC^{m+\theta}({\Bbb R}^3)$ denotes the space
  of all $C^{m}$ functions $u$ on ${\Bbb R}^3$ such that $u$ itself and all
  its partial derivatives of order$\leq m$ are bounded and uniformly
  $\theta$-th order H\"{o}lder continuous in $\Bbb R^3$. We denote $\Pi_1=
  E\circ\Pi$. Then clearly $\Pi_1\in L(C^{m+\theta}(\Bbb S^2),BUC^{m+\theta}
  ({\Bbb R}^3))$, so that there exists a constant $C_0>0$ such that
\begin{equation}
  \|\Pi_1(\rho)\|_{BUC^{m+\theta}({\Bbb R}^3)}
  \le C_0\|\rho\|_{C^{m+\theta}(\Bbb S^2)}
  \quad\mbox{for}\;\;\rho\in C^{m+\theta}(\Bbb S^2).
\end{equation}
  Take a constant $0<\delta<\min\{1/6,1/(3C_0)\}$ and fix it, where $C_0$ is
  the constant in (2.2). We choose a cut-off function $\chi\in C^\infty[0,
  \infty)$ such that
\begin{equation}
  0\le\chi\le1, \quad \chi(\tau)=\left\{
  \begin{array}{ll}
  1,\quad \text{for}\;|\tau|\le\delta,
  \\
  0,\quad \text{for}\;|\tau|\ge3\delta,
  \end{array}
  \right.
  \qquad\text{and}\qquad
  |\chi'(\tau)|\le{2\over3\delta}.
\end{equation}
  We denote
$$
  O_\delta^{m+\theta}(\Bbb S^2)=\{\rho\in
  C^{m+\theta}(\Bbb S^2):||\rho||_{C^{m+\theta}
  ({\Bbb S^2})}<\delta\}.
$$
  Given $\rho\in O_\delta^{m+\theta}(\Bbb S^2)$, we define the {\em
  Hanzawa transformation} $\Phi_\rho:\Bbb R^3\to\Bbb R^3$ as follows:
\begin{equation}
  \Phi_\rho(x)=x+\chi(r-1)\Pi_1(\rho)(x)\frac{x}{r}
  \qquad\mbox{for}\;x\in\Bbb R^3.
\end{equation}
  Using (2.2) and (2.3) we can easily verify that
$$
  \Phi_\rho\in \mbox{Diff}^{m+\theta}(\Bbb R^3,\Bbb R^3)\quad
  \mbox{and}\quad \Phi_\rho(x)=x\quad\mbox{if}\quad\mbox
  {dist}(x,\Bbb S^2)>3\delta.
$$
  We define $\phi_\rho=\Phi_\rho\big|_{\Bbb S^2}$ and $\Gamma_\rho=\mbox{Im}
  (\phi_\rho)$, and denote by $\Omega_\rho$ the domain enclosed by
  $\Gamma_\rho$. Clearly,
$$
  \phi_\rho(\omega)=[1+\rho(\omega)]\omega \quad
  \mbox{for}\;\; \omega\in\Bbb S^2.
$$
  Thus, in the polar coordinates $(r,\omega)$ of $\Bbb R^3$, where $r=|x|$ and
  $\omega=x/|x|$, the hyper-surface $\Gamma_\rho$ has the following equation:
  $r=1+\rho(\omega)$.

  Next, given $\rho\in C([0,T],O_\delta^{m+\theta}(\Bbb S^2)$, for each
  $t\in [0,T]$ we define $\Gamma_{\rho}(t)=\Gamma_{\rho(t)}$ and $\Omega_{\rho}(t)=
  \Omega_{\rho(t)}$. Since our purpose is to study asymptotical stability of
  the radially symmetric stationary solution, later on we always assume the
  initial domain $\Omega_0$ lies in a small neighborhood of $\Omega_s$. More
  precisely, we assume $\Gamma_0:=\partial\Omega_0=\mbox{Im}(\phi_{\rho_0})$
  for some $\rho_0\in O_\delta^{m+\theta}(\Bbb S^2)$.

  Let $\rho$ be as above, and let $\Phi_\rho^i$ be the $i$-th component of
  $\Phi_\rho$, $i=1,2,3$. We denote
$$
  [D\Phi_\rho]_{ij}:=\partial_i\Phi_\rho^j={\partial\Phi_\rho^j\over\partial x_i},
  \quad a_{ij}^\rho(x)=[D\Phi_\rho(x)]^{-1}_{ij}\quad (i,j=1,2,3),
$$
$$
  G_\rho(x)=\det(D\Phi_\rho(x)\big)\quad \mbox{for}\;\; x\in {\Bbb R}^3,
$$
$$
  H_\rho(\omega)=|\phi_\rho|^2\sqrt{1+|\nabla_\omega\phi_\rho|^2}
  \quad \mbox{for}\;\; \omega\in {\Bbb S}^2,
$$
  where $\nabla_\omega$ represents the orthogonal projection of the gradient
  $\nabla_x$ onto the tangent space $T_\omega({\Bbb S}^2)$\footnotemark[1]$^)$.
\footnotetext[1]{$^)$In the coordinate $\omega=\omega(\vartheta,\varphi)=
  (\sin\vartheta\cos\varphi,\sin\vartheta\sin\varphi,\cos\vartheta)$
  ($0\leq\vartheta\leq\pi$, $0\leq\varphi\leq 2\pi$) of the sphere we have
$$
  \nabla_\omega f(\omega)=(\cos\vartheta\cos\varphi,\cos\vartheta
  \sin\varphi,-\sin\vartheta)\partial_\vartheta f(\omega(\vartheta,
  \varphi))+{1\over\sin\vartheta}(-\sin\varphi,\cos\varphi,0)
  \partial_\varphi f(\omega(\vartheta,\varphi)).
$$
  Note also that $\nabla_x f=\displaystyle\frac{\partial f}
  {\partial r}\omega+{1\over r}\nabla_\omega f$.}
  Here and hereafter, for a matrix $A$ we use the notation $A_{ij}$ to denote
  the element of $A$ in the $(i,j)$-th position. From (2.3) we see that
  $[\rho\to\Phi_\rho]\in C^\infty(O^{m+\theta}_\delta({\Bbb S}^2),
  {\rm Diff}^{m+\theta}({\Bbb R}^3,{\Bbb R}^3))$. Thus we have
$$
\left\{
\begin{array}{l}
  [\rho\to a_{ij}^\rho]\in C^\infty(O^{m+\theta}_\delta({\Bbb S}^2),
  C^{m-1+\theta}({\Bbb R}^3)),\quad i,j=1,2,3,\\ [0.2cm]
  [\rho\to G_\rho]\in C^\infty(O^{m+\theta}_\delta({\Bbb S}^2),
  C^{m-1+\theta}({\Bbb R}^3)),\\ [0.2cm]
  [\rho\to H_\rho]\in C^\infty(O^{m+\theta}_\delta({\Bbb S}^2),
  C^{m-1+\theta}({\Bbb S}^2)).
\end{array}
\right.
\eqno{(2.5)}
$$
  We now introduce four partial differential operator
  ${\mathcal A}(\rho)$, $\vec{\mathcal B}(\rho)$, $\vec{\mathcal B}(\rho)
  \cdot$ and $\vec{\mathcal B}(\rho)\otimes$ on ${\Bbb R}^3$ as follows:
$$
  {\mathcal A}(\rho)u(x)=a^\rho_{ij}(x){\partial_j}
  \big(a^\rho_{ik}(x){\partial_k u(x)}\big) \quad
  \mbox{for scalar function}\;\; u,
$$
$$
  \vec{\mathcal B}(\rho)u(x)=\big(a^\rho_{1j}(x)
  {\partial_j u(x)}, a^\rho_{2j}(x){\partial_j u(x)},
  a^\rho_{3j}(x){\partial_j u(x)}\big) \quad
  \mbox{for scalar function}\;\; u,
$$
$$
  \vec{\mathcal B}(\rho)\cdot{\bfv}(x)=a^\rho_{ij}(x)
  {\partial_j v_i(x)} \quad
  \mbox{for vector function}\;\; {\bfv}=(v_1,v_2,v_3).
$$
$$
  \vec{\mathcal B}(\rho)\otimes{\bfv}(x)=(a^\rho_{ik}(x)
  {\partial_k v_j(x)}) \quad
  \mbox{for vector function}\;\; {\bfv}=(v_1,v_2,v_3).
$$
  Here and hereafter we use the convention that repeated indices represent
  summations with respect to these indices, and $\partial_j=\partial/\partial
  x_j$, $j=1,2,3$. These definitions can be respectively briefly rewritten as
  follows:
$$
  {\mathcal A}(\rho)u=(\Delta(u\circ\Phi_\rho^{-1}))\circ\Phi_\rho, \quad
  \vec{\mathcal B}(\rho)u=(\nabla(u\circ\Phi_\rho^{-1}))\circ\Phi_\rho,
$$
$$
  \vec{\mathcal B}(\rho)\cdot{\bfv}=(\nabla\cdot({\bfv}\circ\Phi_\rho^{-1}))
  \circ\Phi_\rho, \quad
  \vec{\mathcal B}(\rho)\otimes{\bfv}=(\nabla\otimes({\bfv}\circ\Phi_\rho^{-1}))
  \circ\Phi_\rho.
$$
  By (2.5) we have
$$
\left\{
\begin{array}{l}
  [\rho\to {\mathcal A}(\rho)]\in C^\infty(O^{m+\theta}_\delta({\Bbb S}^2),
  L(C^{m+\theta}(\overline{\Bbb B}^3),C^{m-2+\theta}(\overline{\Bbb B}^3))),
\\ [0.1cm]
  [\rho\to \vec{\mathcal B}(\rho)]\in
  C^\infty(O^{m+\theta}_\delta({\Bbb S}^2),
  L(C^{m+\theta}(\overline{\Bbb B}^3),
  (C^{m-1+\theta}(\overline{\Bbb B}^3))^3)),
\\ [0.1cm]
  [\rho\to \vec{\mathcal B}(\rho)\cdot]\in
  C^\infty(O^{m+\theta}_\delta({\Bbb S}^2),
  L((C^{m+\theta}(\overline{\Bbb B}^3))^3,
  C^{m-1+\theta}(\overline{\Bbb B}^3))),
\\ [0.1cm]
  [\rho\to \vec{\mathcal B}(\rho)\otimes]\in
  C^\infty(O^{m+\theta}_\delta({\Bbb S}^2),
  L((C^{m+\theta}(\overline{\Bbb B}^3))^3,
  (C^{m-1+\theta}(\overline{\Bbb B}^3))^{3\times 3})).
\end{array}
\right.
\eqno{(2.6)}
$$
  Next we introduce the boundary operator $\vec{\mathcal D}(\rho)$:
  $(C^{m-1+\theta}(\overline{\Bbb B}^3))^3\to C^{m-1+\theta}({\Bbb S}^2)$
$$
  \vec{\mathcal D}(\rho)\bfv=\mbox{\rm tr}_{{\Bbb S}^2}(\bfv)
  \cdot [\omega-{1\over 1\!+\!\rho}\nabla_\omega\rho],
$$
  and the bilinear operator ${\mathcal C}(\rho)$:
  $C^{m+\theta}(\overline{\Bbb B}^3)\times (C^{m-1+\theta}
  (\overline{\Bbb B}^3))^3\to C^{m-1+\theta}(\overline{\Bbb B}^3)$
$$
  {\mathcal C}(\rho)[u,\bfv]=\chi(r-1)\Pi_1(\vec{\mathcal D}(\rho)\bfv)
  \vec{\mathcal B}(\rho)u\cdot\bfe_r.
$$
  Here and hereafter we use the notation $\bfe_r$ to denote the vector function
  on $\Bbb R^3\backslash\{0\}$ defined by $\bfe_r(x)=\omega(x)=x/r$. Note that
  since $\chi(r-1)=0$ for $0\leq r\leq 1-3\delta$, we see that
  ${\mathcal C}(\rho)[u,\bfv]$ is a well-defined function on
  $\overline{\Bbb B}^3$. Again, by (2.5) we have
$$
\left\{
\begin{array}{l}
  [\rho\to \vec{\mathcal D}]\in C^\infty(O^{m+\theta}_\delta({\Bbb S}^2),
  L((C^{m-1+\theta}(\overline{\Bbb B}^3))^3, C^{m-1+\theta}({\Bbb S}^2))),
\\ [0.1cm]
  [\rho\to {\mathcal C}(\rho)]\in C^\infty(O^{m+\theta}_\delta({\Bbb S}^2),
  BL(C^{m+\theta}(\overline{\Bbb B}^3)\times (C^{m-1+\theta}
  (\overline{\Bbb B}^3))^3,C^{m-1+\theta}(\overline{\Bbb B}^3))).
\end{array}
\right.
\eqno{(2.7)}
$$
  Here the notation $BL(\cdot\times\cdot,\cdot)$ denotes the Banach space of
  bounded bilinear mappings with respect to the indicated Banach spaces.

  Let $\bfn$ and $\kappa$ be respectively the unit outward normal and the mean
  curvature of $\Gamma_\rho$ (see (1.5)). We denote
$$
  \widetilde{\bfn}_\rho(x)=\bfn(\phi_\rho(x)), \quad
  \widetilde{\kappa}_\rho(x)=\kappa(\phi_\rho(x)) \quad
  \mbox{for}\;\; x\in {\Bbb S}^2.
$$
  A direct computation shows that
$$
  \widetilde{\bfn}_\rho(x)=\frac{x\cdot[(D\Phi_\rho(x))^{-1}\big]^T}
  {|x\cdot[(D\Phi_\rho(x))^{-1}\big]^T|}=
  {a_{ij}^\rho(x)x_j\bfe_i\over\big|a_{ij}^\rho(x)x_j\bfe_i\big|},
$$
  where
$$
  \bfe_1=(1,0,0),\quad \bfe_2=(0,1,0),\quad \bfe_3=(0,0,1),
$$
  and
$$
  \widetilde{\kappa}_\rho(x)={1\over2}a_{ij}^\rho(x)\partial_j\widetilde
  {\bfn}_\rho^i(x).
$$
  By (2.5) we have
$$
\left\{
\begin{array}{l}
  [\rho\to\widetilde{\kappa}_\rho]\in C^\infty(O^{m+\theta}_\delta({\Bbb S}^2),
  C^{m-2+\theta}({\Bbb S}^2)),\\ [0.2cm]
  [\rho\to\widetilde{\bfn}_\rho]\in C^\infty(O^{m+\theta}_\delta({\Bbb S}^2),
  (C^{m-1+\theta}({\Bbb S}^2))^3).
\end{array}
\right.
\eqno{(2.8)}
$$

  As in \cite{Fried1} we introduce the following vector functions:
$$
  \bfw_1(x)=(0,x_3,-x_2),\quad \bfw_2(x)=(-x_3,0,x_1),\quad
  \bfw_3(x)=(x_2,-x_1,0).
$$
  Then clearly $\bfv\times x=(\bfv\cdot\bfw_1,\bfv\cdot\bfw_2,
  \bfv\cdot\bfw_3)$.

  Let $T$ be a given positive number and consider a function $\rho:
  \;[0,T]\to O_\delta^{m+\theta}({\mathbb S}^2)$. We assume that
  $\rho\in C([0,T],O_\delta^{m+\theta}({\mathbb S}^2))$. Given a
  such $\rho$, we denote
$$
  \Gamma_\rho(t)=\Gamma_{\rho(t)}, \quad
  \Omega_\rho(t)=\Omega_{\rho(t)} \quad (0\leq t\leq T).
$$

  Finally, for $\sigma$, $\bfv$ and $p$ as in (1.1)--(1.9), we denote
$$
  \widetilde{\sigma}=\sigma\circ\Phi_\rho, \quad
  \widetilde{\bfv}=\bfv\circ\Phi_\rho, \quad
  \widetilde{p}=p\circ\Phi_\rho.
$$
  We also denote $\widetilde{\bfw}_j^\rho=\bfw_j\circ\Phi_\rho$, $j=1,2,3$.

  Using these notations, we claim that the Hanzawa transformation transforms
  the equations (1.1)--(1.10) into the following equations, respectively:
$$
  \epsln\partial_t\widetilde{\sigma}-{\mathcal A}(\rho)\widetilde{\sigma}
  -\epsln{\mathcal C}(\rho)[\widetilde{\sigma},\widetilde{\bfv}]
  =-f(\widetilde{\sigma})
  \quad \mbox{in}\;\; {\mathbb B}^3\times \Bbb R_+,
\eqno{(2.9)}
$$
$$
  \vec{\mathcal B}(\rho)\cdot\widetilde{\bfv}=g(\widetilde{\sigma}) \quad
  \mbox{in}\;\; {\mathbb B}^3\times \Bbb R_+,
\eqno{(2.10)}
$$
$$
  -{\mathcal A}(\rho)\widetilde{\bfv}+\vec{\mathcal B}(\rho)\widetilde{p}
  -{1\over 3}\vec{\mathcal B}(\rho)(\vec{\mathcal B}(\rho)\cdot
  \widetilde{\bfv})=0 \quad \mbox{in}\;\; {\mathbb B}^3\times \Bbb R_+,
\eqno{(2.11)}
$$
$$
  \widetilde{\sigma}=1 \quad \mbox{on}\;\; {\mathbb S}^2\times \Bbb R_+,
\eqno{(2.12)}
$$
$$
  \widetilde{\bfT}_\rho(\widetilde{\bfv},\widetilde{p})\widetilde{\bfn}_\rho=
  -\gamma\widetilde{\kappa}_\rho\widetilde{\bfn}_\rho \quad \mbox{on}\;\;
  {\mathbb S}^2\times \Bbb R_+,
\eqno{(2.13)}
$$
$$
  \partial_t\rho=\vec{\mathcal D}(\rho)\widetilde{\bfv} \quad
  \mbox{on}\;\; {\mathbb S}^2\times \Bbb R_+,
\eqno{(2.14)}
$$
$$
  \int_{|x|<1}\widetilde{\bfv}(x)G_\rho(x)dx=0,\quad t>0,
\eqno{(2.15)}
$$
$$
  \int_{|x|<1}\widetilde{\bfv}(x)\cdot\widetilde{\bfw}_j^\rho(x)
  G_\rho(x)dx=0,\quad j=1,2,3, \quad t>0.
\eqno{(2.16)}
$$
$$
  \widetilde{\sigma}(x,0)=\widetilde{\sigma}_0(x) \quad
  \mbox{for}\;\; x\in\Bbb B^3,
\eqno{(2.17)}
$$
$$
  \rho(\omega,0)=\rho_0(\omega) \quad
  \mbox{for}\;\; \omega\in\Bbb S^2.
\eqno{(2.18)}
$$
  Here $\widetilde{\bfT}_\rho(\widetilde{\bfv},\widetilde{p})=
  [\vec{\mathcal B}(\rho)\otimes\widetilde{\bfv}+(\vec{\mathcal B}(\rho)\otimes
  \widetilde{\bfv})^T]-[\widetilde{p}+(2/3)\vec{\mathcal B}(\rho)\cdot
  \widetilde{\bfv}]\bfI$.

  Indeed, it is immediate to see that under the Hanzawa transformation, the
  equations (1.2)--(1.5) and (1.7)--(1.9) are respectively transformed into
  the equations (2.10)--(2.13) and (2.15)--(2.17), and (2.18) is a rewritten
  form of (1.10). In what follows we prove that (1.1) and (1.6) are
  transformations of (2.9) and (2.14).

  Let $\psi_\rho(x,t)=r-1-\rho(\omega,t)$, where $r=|x|$ and $\omega=x/|x|$.
  Then $x\in\Gamma_{\rho}(t)$ if and only if $\psi_\rho(x,t)=0$. It follows
  that the normal velocity of $\Gamma_\rho(t)$ is as follows (see \cite{EscSim}):
$$
  V_n(x,t)=\frac{\partial_t\rho(\omega,t)}{|\nabla_x\psi_\rho(x,t)|} \quad
  \mbox{for}\;\; x\in\Gamma_\rho(t), \quad t>0.
$$
  Moreover, $\bfn(x,t)=\nabla_x\psi_\rho(x,t)/|\nabla_x\psi_\rho(x,t)|$. Hence
  (1.6) can be rewritten as follows:
$$
  \partial_t\rho(\omega,t)=\bfv(x,t)\cdot\nabla_x\psi_\rho(x,t) \quad
  \mbox{for}\;\; x\in\Gamma_\rho(t), \quad t>0,
$$
  where $\omega=x/|x|$. Since $\nabla_x\psi_\rho=\displaystyle
  \frac{\partial\psi_\rho}{\partial r}\omega+{1\over r}\nabla_\omega\psi_\rho$,
  we see that after the Hanzawa transformation, this equation has the
  following form:
$$
  \partial_t\rho(\omega,t)=\widetilde{\bfv}(\omega,t)\cdot
  \Big[\omega-{\nabla_\omega\rho(\omega,t)\over 1\!+\!\rho(\omega,t)}\Big]
  \quad \mbox{for}\;\; \omega\in\Bbb S^2, \quad t>0.
$$
  Recalling the definition of the operator $\vec{\mathcal D}(\rho)$, we see that
  the equation (2.14) follows. Next, by differentiating the relation
  $\widetilde{\sigma}=\sigma\circ\Phi_\rho$ in $t$ and using the equations (1.1),
  (2.4) and (2.14) we see that
$$
\begin{array}{rcl}
  \partial_t\widetilde{\sigma}&=&\partial_t\sigma\circ\Phi_\rho+
  \partial_t\Phi_\rho\cdot(\nabla_x\sigma\circ\Phi_\rho) \\ [0.2cm]
  &=&
\displaystyle {1\over\epsln}\Delta\sigma\circ\Phi_\rho-
  {1\over\epsln}f(\sigma\circ\Phi_\rho)+\chi(r-1)\Pi_1(\partial_t\rho)
  \bfe_r\cdot(\nabla_x\sigma\circ\Phi_\rho) \\ [0.2cm]
  &=&
\displaystyle {1\over\epsln}{\mathcal A}(\rho)\widetilde{\sigma}-
  {1\over\epsln}f(\widetilde{\sigma})+\chi(r-1)
  \Pi_1(\vec{\mathcal D}(\rho)\widetilde{\bfv})
  \bfe_r\cdot{\mathcal B}(\rho)\widetilde{\sigma} \\ [0.2cm]
  &=&
\displaystyle {1\over\epsln}{\mathcal A}(\rho)\widetilde{\sigma}-
  {1\over\epsln}f(\widetilde{\sigma})+
  {\mathcal C}(\rho)[\widetilde{\sigma},\widetilde{\bfv})].
\end{array}
$$
  Hence (2.9) follows.

  The above deduction yields the following lemma:
\medskip

  {\bf Lemma 2.1}\ \ {\em If $(\sigma,\bfv,p,\rho)$ is a solution of the problem
  $(1.1)$--$(1.10)$, then by letting $\widetilde{\sigma}=\sigma\circ\Phi_\rho$,
  $\widetilde{\bfv}=\bfv\circ\Phi_\rho$ and $\widetilde{p}=p\circ\Phi_\rho$, we
  have that $(\widetilde{\sigma},\widetilde{\bfv},\widetilde{p},\rho)$ is a
  solution of the problem $(2.9)$--$(2.18)$. Conversely, If $(\widetilde{\sigma},
  \widetilde{\bfv},\widetilde{p},\rho)$ is a solution of the problem
  $(2.9)$--$(2.18)$, then by letting $\sigma=\widetilde{\sigma}\circ
  \Phi_\rho^{-1}$, $\bfv=\widetilde{\bfv}\circ\Phi_\rho^{-1}$ and $p=
  \widetilde{p}\circ\Phi_\rho^{-1}$, we have that  $(\sigma,\bfv,p,\rho)$ is a
  solution of the problem $(1.1)$--$(1.10)$.}
\medskip

  {\em Proof}:\ \ The above deduction shows that if $(\sigma,\bfv,p,\rho)$ is a
  solution of $(1.1)$--$(1.10)$, then $(\widetilde{\sigma},\widetilde{\bfv},
  \widetilde{p},\rho)$ satisfies $(2.9)$--$(2.18)$. The converse can be
  similarly verified. $\qquad$ $\Box$
\medskip

  We now proceed to reduce the problem (2.9)--(2.18) into evolution equations
  only in $\widetilde{\sigma}$ and $\rho$. The idea is to solve equations
  (2.10), (2.11), (2.13), (2.15) and (2.16) to get $\widetilde{\bfv}$ and
  $\widetilde{p}$ as functionals of $\widetilde{\sigma}$ and $\rho$, and next
  substitute $\widetilde{\bfv}$ obtained in this way into equations (2.9) and
  (2.14). Thus, for given $\rho\in O^{m+\theta}_\delta(\Bbb S^2)$ we consider
  the following boundary value problem:
$$
  \vec{\mathcal B}(\rho)\cdot\widetilde{\bfv}=\varphi \quad
  \mbox{in}\;\; \Bbb B^3,
\eqno{(2.19)}
$$
$$
  -{\mathcal A}(\rho)\widetilde{\bfv}+\vec{\mathcal B}(\rho)\widetilde{p}
  =\bfg \quad \mbox{in}\;\; \Bbb B^3,
\eqno{(2.20)}
$$
$$
  \widetilde{\bfT}_\rho(\widetilde{\bfv},\widetilde{p})\widetilde{\bfn}_\rho=
  \bfh \quad \mbox{on}\;\; \Bbb S^2,
\eqno{(2.21)}
$$
$$
  \int_{|x|<1}\widetilde{\bfv}(x)G_\rho(x)dx=0,
\eqno{(2.22)}
$$
$$
  \int_{|x|<1}\widetilde{\bfv}(x)\cdot\widetilde{\bfw}_j^\rho(x)
  G_\rho(x)dx=0,\quad j=1,2,3,
\eqno{(2.23)}
$$
  where $\varphi\in C^{m\!-\!k\!-\!1\!+\!\theta}(\overline{\Bbb B}^3)$, $\bfg
  \in (C^{m\!-\!k\!-\!2\!+\!\theta}(\overline{\Bbb B}^3))^3$ and $\bfh\in
  (C^{m\!-\!k\!-\!1\!+\!\theta}({\Bbb S}^2))^3$ for some $0\!\leq k\leq\!
  m\!-\!2$.
\medskip

  {\bf Lemma 2.2}\ \ {\em Let $\delta$ be sufficiently small and let $\rho\in
  O^{m+\theta}_\delta(\Bbb S^2)$ be given. A necessary and sufficient condition
  for $(2.19)$--$(2.23)$ to have a solution is that $\varphi$, $\bfg$ and $\bfh$
  satisfy the following relations:
$$
  \int_{|x|<1}\big(\bfg(x)-{1\over 3}\vec{\mathcal B}(\rho)\varphi(x)\big)
  \cdot\widetilde{\bfw}_j^\rho(x) G_\rho(x)dx+
  \int_{|x|=1}\bfh(x)\cdot\widetilde{\bfw}_j^\rho(x)
  H_\rho(x) dS_x=0, \quad j=1,2,3,
\eqno{(2.24)}
$$
$$
  \int_{|x|<1}\big(\bfg(x)-{1\over 3}\vec{\mathcal B}(\rho)\varphi(x)\big)
  \cdot{\bfe}_j G_\rho(x)dx+
  \int_{|x|=1}\bfh(x)\cdot{\bfe}_j
  H_\rho(x) dS_x=0, \quad j=1,2,3.
\eqno{(2.25)}
$$
  If this condition is satisfied, then $(2.19)$--$(2.23)$ has a unique solution
  $(\widetilde{\bfv},\widetilde{p})\in (C^{m-k+\theta}(\overline{\Bbb B}^3))^3
  \times C^{m-k-1+\theta}(\overline{\Bbb B}^3)$. Moreover, we have $\widetilde
  {\bfv}=\vec{\mathcal P}(\rho)\varphi+{\bf Q}(\rho)\bfg+{\bf R}(\rho)\bfh$, where
$$
\left\{
\begin{array}{l}
\displaystyle
  \vec{\mathcal P}\in\bigcap_{k=0}^{m-2} C^\infty(O^{m+\theta}_\delta(\Bbb S^2),
  L(C^{m-k-1+\theta}(\overline{\Bbb B}^3),
  (C^{m-k+\theta}(\overline{\Bbb B}^3))^3),\\ [0.3cm]
\displaystyle
  {\bf Q}\in\bigcap_{k=0}^{m-2} C^\infty(O^{m+\theta}_\delta(\Bbb S^2),
  L((C^{m-k-2+\theta}(\overline{\Bbb B}^3))^3,
  (C^{m-k+\theta}(\overline{\Bbb B}^3))^3)),\\ [0.3cm]
\displaystyle
  {\bf R}\in\bigcap_{k=0}^{m-2} C^\infty(O^{m+\theta}_\delta(\Bbb S^2),
  L((C^{m-k-1+\theta}({\Bbb S}^2))^3,
  (C^{m-k+\theta}(\overline{\Bbb B}^3))^3)).
\end{array}
\right.
\eqno{(2.26)}
$$
}

  {\em Proof}:\ \ See Lemmas 2.3 and 2.4 of \cite{WuCui}. $\qquad$ $\Box$
\medskip

  Now, let $\widetilde{\sigma}\in C^{m+\theta}(\overline{\Bbb B}^3)$, $\rho\in
  O^{m+\theta}_\delta(\Bbb S^2)$ and we consider the system of equations
  (2.10), (2.11), (2.13), (2.15) and (2.16). These equations can be rewritten
  in the form of (2.19)--(2.23), with
$$
  \varphi=g(\widetilde{\sigma}), \quad
  \bfg={1\over 3}\vec{\mathcal B}(\rho)g(\widetilde{\sigma}), \quad
  \bfh=-\gamma\widetilde{\kappa}_\rho\widetilde{\bfn}_\rho.
\eqno{(2.27)}
$$
  As was shown in \cite{WuCui}, the relations (2.24) and (2.25) are satisfied
  by these functions. Besides, it is obvious that $\varphi\in C^{m+\theta}
  (\overline{\Bbb B}^3)\subseteq C^{m-2+\theta}(\overline{\Bbb B}^3)$ and
  $\bfg\in (C^{m-1+\theta}(\overline{\Bbb B}^3))^3\subseteq (C^{m-3+\theta}
  (\overline{\Bbb B}^3))^3$. Furthermore, by (2.8) we see that $\bfh\in
  (C^{m-2+\theta}({\Bbb S}^2))^3$. Hence, by Lemma 2.2 (with $k=1$) we infer
  that these equations have a unique solution $(\widetilde{\bfv},\widetilde{p})
  \in (C^{m-1+\theta}(\overline{\Bbb B}^3))^3\times C^{m-2+\theta}
  (\overline{\Bbb B}^3)$, and
$$
  \widetilde{\bfv}=\vec{\mathcal V}(\widetilde{\sigma},\rho)\equiv
  \vec{\mathcal P}(\rho)g(\widetilde{\sigma})
  +{1\over 3}{\mathbf Q}(\rho)\vec{\mathcal B}(\rho)g(\widetilde{\sigma})
  -\gamma {\mathbf R}(\rho)({\mathcal K}(\rho)\vec{\mathcal N}(\rho)).
\eqno{(2.28)}
$$
  where ${\mathcal K}(\rho)=\widetilde{\kappa}_\rho$ and $\vec{\mathcal N}
  (\rho)=\widetilde{\bfn}_\rho$. We note that
$$
  {\mathcal K}\in C^\infty(O^{m+\theta}_\delta({\Bbb S}^2),
  C^{m-2+\theta}({\Bbb S}^2)), \quad
  \vec{\mathcal N}\in C^\infty(O^{m+\theta}_\delta({\Bbb S}^2),
  (C^{m-1+\theta}({\Bbb S}^2))^3).
\eqno{(2.29)}
$$
  Substituting the above expression of $\widetilde{\bfv}$ into (2.9) and (2.14),
  and introducing operators ${\mathcal F}:C^{m+\theta}(\overline{\Bbb B}^3)
  \times O^{m+\theta}_\delta({\Bbb S}^2)\to C^{m-2+\theta}(\overline{\Bbb B}^3)$
  and ${\mathcal G}:C^{m+\theta}(\overline{\Bbb B}^3)\times
  O^{m+\theta}_\delta({\Bbb S}^2)\to C^{m-1+\theta}({\Bbb S}^2)$ respectively by
$$
  {\mathcal F}(\widetilde{\sigma},\rho)={1\over\epsln}{\mathcal A}(\rho)
  \widetilde{\sigma}-\chi(|x|-1)\Pi_1\big(\vec{\mathcal D}(\rho)
  \vec{\mathcal V}(\widetilde{\sigma},\rho)\big)
  \vec{\mathcal B}(\rho)\widetilde{\sigma}\cdot\bfe_r
  -{1\over\epsln}f(\widetilde{\sigma}),
\eqno{(2.30)}
$$
$$
  {\mathcal G}(\widetilde{\sigma},\rho)=\vec{\mathcal D}(\rho)
  \vec{\mathcal V}(\widetilde{\sigma},\rho)=
  {\rm tr}_{{\Bbb S}^2}
  \big[\vec{\mathcal V}(\widetilde{\sigma},\rho)\big]\cdot
  \Big[\omega-{1\over 1+\rho}\nabla_\omega\rho\Big],
\eqno{(2.31)}
$$
  (for $\widetilde{\sigma}\in C^{m+\theta}(\overline{\Bbb B}^3)$ and $\rho\in
  O^{m+\theta}_\delta({\Bbb S}^2)$), where as before $\omega$ represents the
  variable in ${\Bbb S}^2$ and $\omega(x)=x/|x|$ for $x\in \overline{\Bbb B}^3
  \backslash\{0\}$, we see that the problem (2.9)--(2.18) is reduced into the
  following problem:
$$
\left\{
\begin{array}{l}
   \partial_t\widetilde{\sigma}={\mathcal F}(\widetilde{\sigma},\rho) \quad
   \mbox{in}\;\;{\Bbb B}^3\times {\Bbb R}_+,\\ [0.1cm]
   \partial_t\rho={\mathcal G}(\widetilde{\sigma},\rho) \quad
   \mbox{on}\;\;{\Bbb S}^2\times {\Bbb R}_+,\\ [0.1cm]
   \widetilde{\sigma}=1 \quad \mbox{on}\;\; {\mathbb S}^2\times \Bbb R_+,
   \\ [0.1cm]
   \widetilde{\sigma}|_{t=0}=\widetilde{\sigma}_0 \quad
   \mbox{in}\;\; {\Bbb B}^3,\\ [0.1cm]
   \rho|_{t=0}=\rho_0\quad
   \mbox{on}\;\; {\Bbb S}^2.
\end{array}
\right.
\eqno{(2.32)}
$$

  We summarize:
\medskip

  {\bf Lemma 2.3}\ \ {\em Let $(\widetilde{\sigma},\widetilde{\bfv},
  \widetilde{p},\rho)$ be a solution of the problem $(2.9)$--$(2.18)$. Then
  $(\widetilde{\sigma},\rho)$ is a solution of the problem $(2.32)$.
  Conversely, if $(\widetilde{\sigma},\rho)$ is a solution of $(2.32)$, then
  by letting $(\widetilde{\bfv},\widetilde{p})$ be the unique solution of the
  problem $(2.19)$--$(2.23)$ in which $\varphi$, $\bfg$ and $\bfh$ are given
  by $(2.27)$, we have that $(\widetilde{\sigma},\widetilde{\bfv},\widetilde{p},
  \rho)$ is a solution of $(2.9)$--$(2.18)$.} $\qquad$ $\Box$
\medskip

  The problem (2.30) can be rewritten as an initial value problem of a
  differential equation in a Banach space. For this purpose we denote
$$
  {\Bbb X}=C^{m-2+\theta}(\overline{\Bbb B}^3)\times C^{m-1+\theta}({\Bbb S}^2),
  \qquad
  {\Bbb X}_0=(C^{m+\theta}(\overline{\Bbb B}^3)\cap C_0(\overline{\Bbb B}^3))
  \times C^{m+\theta}({\Bbb S}^2),
$$
$$
  {\Bbb O}_\delta=(C^{m+\theta}(\overline{\Bbb B}^3)\cap C_0(\overline{\Bbb B}^3))
  \times O^{m+\theta}_\delta({\Bbb S}^2),
$$
  where $C_0(\overline{\Bbb B}^3)=\{u\in C(\overline{\Bbb B}^3):u|_{{\Bbb S}^2}
  =0\}$, and define a bounded nonlinear operator ${\Bbb F}$ in ${\Bbb X}$
  with domain ${\Bbb O}_\delta$ (i.e., ${\Bbb F}:{\Bbb O}_\delta\to {\Bbb X}$) as
  follows:
$$
  {\Bbb F}(U)=({\mathcal F}(u+1,\rho),{\mathcal G}(u+1,\rho)) \quad
  \mbox{for}\;\; U=(u,\rho)\in {\Bbb O}_\delta.
\eqno{(2.33)}
$$
  Then (2.32) can be rewritten as an initial value problem of a differential
  equation in ${\Bbb X}$:
$$
\left\{
\begin{array}{l}
   \displaystyle{dU\over dt}={\Bbb F}(U) \quad
   \mbox{for}\;\;t>0,\\ [0.2cm]
   U|_{t=0}=U_0,
\end{array}
\right.
\eqno{(2.34)}
$$
  where $U_0=(\widetilde{\sigma}_0-1,\rho_0)$. The relation between solutions
  of (2.32) and (2.34) is that $U=(\widetilde{\sigma}-1,\rho)$.

  We note that ${\Bbb X}_0$ is not dense in ${\Bbb X}$.
\medskip

\section{Linearization of ${\Bbb F}(U)$}
\setcounter{equation}{0}

  Let $(\sigma_s,\bfv_s,p_s,R_s)$ be the radially symmetric stationary solution
  of the problem (1.1)--(1.10) (recall that $R_s=1$) and denote $U_s=(\sigma_s
  \!-\!1,0)$. Then $U_s$ is a stationary solution of the differential equation
  in (2.31), so that ${\Bbb F}(U_s)=0$.

  From (2.6), (2.7) and (2.29) it can be easily seen that ${\Bbb F}\in
  C^\infty({\Bbb O},{\Bbb X})$, where ${\Bbb O}$ is regarded as an open subset
  of ${\Bbb X}_0$. It follows that the Fr\'echet derivative $D{\Bbb F}\in
  C^\infty({\Bbb O},L({\Bbb X}_0,{\Bbb X}))$.
  In this section we first derive a useful expression of $D{\Bbb F}(U_s)$, and
  next use it to prove that $D{\Bbb F}(U_s)$ is an infinitesimal generator of an
  analytic semigroup in ${\Bbb X}$ with domain ${\Bbb X}_0$.

  By (2.33) we see that for $V=(v,\eta)\in \Bbb X_0$, we have
\begin{equation}
  D\Bbb F(U_s)V=(D_{\widetilde\sigma}\mathcal F(\sigma_s,0)v
  +D_\rho\mathcal F(\sigma_s,0)\eta,
  D_{\widetilde\sigma}\mathcal G(\sigma_s,0)v
  +D_\rho\mathcal G(\sigma_s,0)\eta),
\end{equation}
  where $D_{\widetilde\sigma}\mathcal F$ and $D_\rho\mathcal F$ represent Fr\'echet
  derivatives of $\mathcal F(\widetilde\sigma,\rho)$ in $\widetilde\sigma$ and $\rho$,
  respectively, and similarly for $D_{\widetilde\sigma}\mathcal G$ and $D_\rho\mathcal G$.
  In what follows we deduce expressions of these Fr\'echet derivatives.

  We first note that, clearly,
\begin{equation}
  \mathcal A(0)u=\Delta u,\qquad  \vec{\mathcal B}(0)u=\nabla u,\qquad
  \vec{\mathcal B}(0)\cdot\bfv=\nabla\cdot\bfv,\qquad
  \vec{\mathcal B}(0)\otimes\bfv=\nabla\otimes\bfv,
\end{equation}
\begin{equation}
  \vec{\mathcal D}(0)\bfv={\rm tr}_{\Bbb S^2}(\bfv)\cdot\bfn_0,\qquad
  \mathcal K(0)=\widetilde{\kappa}_\rho|_{\rho=0}=1,\qquad
  \vec{\mathcal N}(0)=\widetilde{\bfn}_\rho|_{\rho=0}=\bfn_0,
\end{equation}
\begin{equation}
  G_\rho(x)|_{\rho=0}=1,\qquad  H_\rho(x)|_{\rho=0}=1,
  \qquad \widetilde\bfw_j^\rho(x)|_{\rho=0}=\bfw_j(x).
\end{equation}
  In (3.3) $\bfn_0$ denotes the unit outward normal of the unit sphere
  $\Bbb S^2$, and this notation will be used throughout the remaining part of
  this paper. We also denote
$$
  M(\eta)=\chi(r-1)\Pi_1(\eta), \qquad
  U(\eta)=\lim_{\epsilon\to0}\frac{G_{\epsilon\eta}-1}{\epsilon}, \qquad
  \bfW_j(\eta)=\lim_{\epsilon\to0}
  \frac{\widetilde\bfw_j^{\epsilon\eta}-\bfw_j}{\epsilon}.
$$
  They are evidently linear operators in $\eta$.
\medskip

  {\bf Lemma 3.1}\ \ {\em We have
\begin{equation}
  [\mathcal A'(0)\eta]\sigma_s=-[\Delta-f'(\sigma_s(r))]
  \big[\sigma_s'(r)M(\eta)\big]\quad \mbox{\rm in}\;\;{\Bbb B^3},
\end{equation}
\begin{equation}
  [\vec{\mathcal B}'(0)\eta]\cdot\bfv_s=-\nabla\cdot[v_s'(r)M(\eta)\bfe_r]
  +g'(\sigma_s)\sigma_s'(r)M(\eta)\quad \mbox{\rm in}\;\;{\Bbb B^3},
\end{equation}
\begin{eqnarray}
  -[\mathcal A'(0)\eta]\bfv_s+[\vec{\mathcal B}'(0)\eta]p_s
  &-&{1\over3}[\vec{\mathcal B}'(0)\eta]g(\sigma_s)
  ={1\over3}\nabla[g'(\sigma_s)\sigma_s'(r)M(\eta)]
\nonumber \\
  &&-\nabla[p_s'(r)M(\eta)]+\Delta[v_s'(r)M(\eta)\bfe_r]
  \quad \mbox{\rm in}\;\;{\Bbb B^3},
\end{eqnarray}
\begin{eqnarray}
  \Big[[\vec{\mathcal B}'(0)\eta]\otimes\bfv_s+\Big([\vec{\mathcal B}'(0)\eta]
  \otimes\bfv_s\Big)^T\Big]\bfn_0 &=& -\Big[\nabla\otimes
  [v_s'(r)M(\eta)\bfe_r]+\Big(\nabla\otimes
  [v_s'(r)M(\eta)\bfe_r]\Big)^T\Big]\bfn_0
\nonumber\\
  &&+\Big[p_s'(r)+{2\over3}g'(\sigma_s)\sigma_s'(r)-4g(1)\Big]\eta\,\bfn_0
  \quad \mbox{\rm on}\;\;{\Bbb S^2},
\end{eqnarray}
\begin{equation}
  [\vec{\mathcal D}'(0)\eta]\bfv=-{\rm tr}_{\Bbb S^2}(\bfv)\cdot
  \nabla_\omega\eta\quad \mbox{\rm in}\;\;{\Bbb B^3},
\end{equation}
\begin{equation}
  \int_{|x|<1}\big[v_s'(r)M(\eta)\bfe_r+U(\eta)\bfv_s\big]\,dx=0,
\end{equation}
\begin{equation}
  \int_{|x|<1}\bfv_s\cdot\bfW_j(\eta)\,dx=0,
  \quad j=1,2,3.
\end{equation}
}

\medskip
  {\em Proof}:\ \ (3.5) follows from (5.8) of \cite{Cui2}. To prove (3.6)
  we denote $\sigma_{s,\epsilon\eta}=\sigma_s\circ\Phi_{\epsilon\eta}$ and
  $\bfv_{s,\epsilon\eta}=\bfv_s\circ\Phi_{\epsilon\eta}$. By making
  Hanzawa transformation to the equation $\nabla\cdot\bfv_{s}=g(\sigma_{s})$
  we have
\begin{equation}
  \vec{\mathcal B}(\epsilon\eta)\cdot
  \bfv_{s,\epsilon\eta}=g(\sigma_{s,\epsilon\eta}).
\end{equation}
  Since $\bfv_{s,\epsilon\eta}|_{\epsilon=0}=\bfv_{s}$ and
  $\sigma_{s,\epsilon\eta}|_{\epsilon=0}=\sigma_{s}$, we get
$$
   [\vec{\mathcal B}(\epsilon\eta)-\vec{\mathcal B}(0)]\cdot
  \bfv_{s,\epsilon\eta}+\vec{\mathcal B}(0)\cdot[\bfv_{s,\epsilon\eta}
  -\bfv_s]=g(\sigma_{s,\epsilon\eta})-g(\sigma_s).
$$
  Dividing both sides with $\epsilon$, then letting $\epsilon\to 0$ and
  using the relations
\begin{equation}
  \lim_{\epsilon\to0}{\sigma_{s,\epsilon\eta}-\sigma_s\over\epsilon}=
  \sigma_s'(r)M(\eta), \qquad
  \lim_{\epsilon\to0}{\bfv_{s,\epsilon\eta}-\bfv_s
  \over\epsilon}=v_s'(r)M(\eta)\bfe_r,
\end{equation}
  we see that (3.6) follows. To prove (3.7) we denote $p_{s,\epsilon\eta}=
  p_s\circ\Phi_{\epsilon\eta}$. By making Hanzawa transformation to the equation
  $-\Delta\bfv_s+\nabla p_s-{1\over3}\nabla(g(\sigma_s))=0$ we get
\begin{equation}
  -\mathcal A(\epsilon\eta)\bfv_{s,\epsilon\eta}+\vec{\mathcal B}
  (\epsilon\eta)p_{s,\epsilon\eta}-{1\over3}\vec{\mathcal B}
  (\epsilon\eta)(g(\sigma_{s,\epsilon\eta}))=0.
\end{equation}
  Since $p_{s,\epsilon\eta}|_{\epsilon=0}=p_{s}$ and
\begin{equation}
  \lim_{\epsilon\to0}{p_{s,\epsilon\eta}-p_s\over\epsilon}=p_s'(r)M(\eta),
\end{equation}
  by a similar argument as before we obtain (3.7).

  Next we prove (3.8). We denote $\widetilde\bfe_r^{\epsilon\eta}=\bfe_r\circ
  \Phi_{\epsilon\eta}$. By making Hanzawa transformation to the equation
  $\bfT(\bfv_s,p_s)\big|_{\Bbb S^2}\bfn_0=-\gamma\bfn_0$ and noticing that
  $\bfe_r\big|_{\Bbb S^2}=\bfn_0$, we get
\begin{equation}
  \widetilde{\bfT}_{\epsilon\eta}(\bfv_{s,\epsilon\eta},p_{s,\epsilon\eta})
  \widetilde\bfe_r^{\epsilon\eta}+\gamma\widetilde\bfe_r^{\epsilon\eta}=0 \qquad
  \mbox{on}\;\;\Phi_{\epsilon\eta}^{-1}(\Bbb S^2).
\end{equation}
  By (3.13) and (3.15) we have $\bfv_{s,\epsilon\eta}=\bfv_{s}+\epsilon
  v_s'(r)M(\eta)\bfe_r+o(\epsilon)$ and $p_{s,\epsilon\eta}=p_{s}+\epsilon
  p_s'(r)M(\eta)+o(\epsilon)$. Thus
\begin{eqnarray*}
  \widetilde{\bfT}_{\epsilon\eta}(\bfv_{s,\epsilon\eta},p_{s,\epsilon\eta})
  &=&\Big[\vec{\mathcal B}(\epsilon\eta)\otimes\bfv_{s,\epsilon\eta}
  +\Big(\vec{\mathcal B}(\epsilon\eta)\otimes\bfv_{s,\epsilon\eta}\Big)^T
  \Big]-\Big[p_{s,\epsilon\eta}+{2\over3}\vec{\mathcal B}(\epsilon\eta)\cdot
  \bfv_{s,\epsilon\eta}\Big]\bfI
  \nonumber\\
  &=&\Big[\vec{\mathcal B}(0)\otimes\bfv_{s}
  +\Big(\vec{\mathcal B}(0)\otimes\bfv_{s}\Big)^T\Big]
  -\Big[p_{s}+{2\over3}\vec{\mathcal B}(0)\cdot\bfv_{s}\Big]\bfI
  \nonumber\\
  &&+\epsilon\Big\{[\vec{\mathcal B}'(0)\eta]\otimes\bfv_{s}
  +\Big([\vec{\mathcal B}'(0)\eta]\otimes\bfv_{s}\Big)^T
  -{2\over3}[\vec{\mathcal B}'(0)\eta]\cdot\bfv_{s}\bfI
  \nonumber\\
  &&+\vec{\mathcal B}(0)\otimes [v_s'(r)M(\eta)\bfe_r]
  +\Big(\vec{\mathcal B}(0)\otimes [v_s'(r)M(\eta)\bfe_r]\Big)^T
  \nonumber\\
  &&-p_s'(r)M(\eta)\bfI-{2\over3}\vec{\mathcal B}(0)\cdot
  [v_s'(r)M(\eta)\omega]\bfI\Big\}+o(\epsilon).
\nonumber
\end{eqnarray*}
  Noticing that
$$
  \big[\vec{\mathcal B}(0)\otimes\bfv_{s}+\big(\vec{\mathcal B}(0)\otimes
  \bfv_{s}\big)^T\big]-\big[p_{s}+{2\over3}\vec{\mathcal B}(0)\cdot
  \bfv_{s}\big]\bfI
  =\bfT(\bfv_s,p_s),
$$
  and denoting by $\bfL(\eta)$ the expression in the braces, we see that the
  above result can be briefly rewritten as follows:
\begin{equation}
  \widetilde{\bfT}_{\epsilon\eta}(\bfv_{s,\epsilon\eta},p_{s,\epsilon\eta})
  =\bfT(\bfv_s,p_s)+\epsilon\bfL(\eta)+o(\epsilon).
\end{equation}
  Since for $x\in\Phi_{\epsilon\eta}^{-1}(\Bbb S^2)$ we have
$$
  \Phi_{\epsilon\eta}(x)=x+\epsilon\chi(r-1)\Pi_1(\eta)(x)x/r=
  [r+\epsilon\chi(r-1)\Pi_1(\eta)(x)]\omega(x),
$$
  where $\omega(x)=x/r$, and $\Phi_{\epsilon\eta}(x)\in\Bbb S^2$, we see that
  $\Phi_{\epsilon\eta}(x)=\omega(x)$ for all $x\in\Phi_{\epsilon\eta}^{-1}
  (\Bbb S^2)$. This implies that $\widetilde\bfe_r^{\epsilon\eta}(x)=
  \bfe_r(\omega(x))=\bfn_0(\omega(x))$ for all $x\in\Phi_{\epsilon\eta}^{-1}
  (\Bbb S^2)$. Hence, from (3.16) and (3.17) we get
\begin{eqnarray}
  0&\;=\;&[\widetilde{\bfT}_{\epsilon\eta}
  (\bfv_{s,\epsilon\eta},p_{s,\epsilon\eta})\widetilde\bfe_r^{\epsilon\eta}
  +\gamma\widetilde\bfe_r^{\epsilon\eta}]\big|_{\Phi^{-1}_{\epsilon\eta}(\Bbb S^2)}
  \nonumber\\
  &\;=\;&[\bfT(\bfv_s,p_s)\bfn_0\circ\omega+\gamma\bfn_0\circ\omega]
  \big|_{\Phi^{-1}_{\epsilon\eta}(\Bbb S^2)}+\epsilon\,\bfL(\eta)
  \big|_{\Phi^{-1}_{\epsilon\eta}(\Bbb S^2)}+o(\epsilon)
  \nonumber\\
  &\;=\;&\Big[\bfT(\bfv_s,p_s)\bfn_0+\gamma\bfn_0-\epsilon\,
  \eta{\partial\over\partial r}\Big(\bfT(\bfv_s,p_s)\bfn_0\Big)\Big]
  \Big|_{\Bbb S^2}+\epsilon\,\bfL(\eta)\Big|_{\Bbb S^2}+o(\epsilon).
\end{eqnarray}
  Points on $\Phi^{-1}_{\epsilon\eta}(\Bbb S^2)$ and $\Bbb S^2$ such that
  the last equality holds are related by the relation $\omega=\omega(x)$ for
  $x\in\Phi^{-1}_{\epsilon\eta}(\Bbb S^2)$ and $\omega\in\Bbb S^2$, and in
  getting the last equality we used the following relations:
$$
  [\bfT(\bfv_s,p_s)\bfn_0\circ\omega+\gamma\bfn_0\circ\omega]
  \big|_{\Phi^{-1}_{\epsilon\eta}(\Bbb S^2)}
  =\Big[\bfT(\bfv_s,p_s)\bfn_0+\gamma\bfn_0-\epsilon\,
  \eta{\partial\over\partial r}\Big(\bfT(\bfv_s,p_s)\bfn_0\Big)\Big]
  \Big|_{\Bbb S^2}+o(\epsilon),
$$
$$
  \bfL(\eta)\big|_{\Phi^{-1}_{\epsilon\eta}(\Bbb S^2)}
  =\bfL(\eta)\big|_{\Bbb S^2}+O(\epsilon).
$$
  The proof of the first relation uses a similar argument as that used in
  (4.29) of \cite{FriedHu1}, and the second relation is immediate. Since
  $[\bfT(\bfv_s,p_s)\bfn_0+\gamma\bfn_0]|_{{\Bbb S}^2}=0$,
  $M(\eta)|_{{\Bbb S}^2}=\eta$, and by the result in Appendix A of
  \cite{WuCui} we have
$$
  {\partial\over\partial r}\Big(\bfT(\bfv_s,p_s)\bfn_0\Big)\Big|_{\Bbb S^2}
  =[2v_s''(r)-p_s'(r)-{2\over3}g'(\sigma_s)\sigma_s'(r)]\Big|_{r=1}
  \bfn_0=-4g(1)\bfn_0,
$$
  by dividing (3.18) with $\epsilon$, then letting $\epsilon\to 0$ and using
  (3.6), we see that (3.8) follows.

  Finally, (3.9) is immediate, and (3.10), (3.11) follow from the relations
  $\displaystyle\int_{\Bbb S^2}\bfv_s\,dx=0$, $\displaystyle\int_{\Bbb S^2}
  \bfv_s\times{x}\,d\,x=0$ and a similar argument as above, which we omit here.
  This completes the proof of Lemma 3.1. $\qquad\Box$
\medskip

  {\bf Lemma 3.2}\ \ {\em For $v\in C^{m+\theta}(\overline{\Bbb B}^3)$ and
  $\eta\in C^{m+\theta}(\Bbb S^2)$ we have
\begin{eqnarray}
  D_{\tilde{\sigma}}\vec{\mathcal V}(\sigma_s,0)v\;&=&
  \vec{\mathcal P}(0)[g'(\sigma_s)v]+{1\over 3}{\mathbf Q}(0)
  \{\nabla[g'(\sigma_s)v]\},
  \\
  D_\rho\vec{\mathcal V}(\sigma_s,0)\eta\;&=&\;v_s'(r)M(\eta)\bfe_r-
  \vec{\mathcal P}(0)\big[g'(\sigma_s)\sigma_s'(r)M(\eta)\big]
  -{1\over3}\mathbf Q(0)\big\{\nabla[g'(\sigma_s)\sigma_s'(r)M(\eta)
  ]\big\}
  \nonumber\\
  &&+{\mathbf R}(0)\big[\gamma(\eta+{1\over2}
  \Delta_\omega\eta)\bfn_0-2g(1)\nabla_\omega\eta+4g(1)\eta
  \bfn_0\big].
\end{eqnarray}
}
\medskip

  {\em Proof:}\ \ From the definition of $\vec{\mathcal V}(\widetilde\sigma,\rho)$
  it is clear that
$$
  D_{\widetilde\sigma}\vec{\mathcal V}(\sigma_s,0)v
  =\vec{\mathcal P}(0)[g'(\sigma_s)v]+{1\over 3}{\mathbf Q}(0)
  \vec{\mathcal B}(0)[g'(\sigma_s)v].
$$
  Since $\vec{\mathcal B}(0)=\nabla$, we see that (3.19) follows.

  To compute $\bfV\equiv D_\rho\vec{\mathcal V}(\sigma_s,0)\eta$ we denote
  $\widetilde{\bfv}=\vec{\mathcal V}(\sigma_s,\epsilon\eta)$, where $\eta\in
  C^{m+\theta}(\Bbb S^2)$ is given. By the definition of $\vec{\mathcal V}
  (\widetilde\sigma,\rho)$ we see that there exists a function $\widetilde{p}
  \in C^{m-1+\theta}(\overline{\Bbb B}^3)$ such that $(\widetilde{\bfv},
  \widetilde{p})$ is the unique solution of the problem
\begin{equation}
  \vec{\mathcal B}(\epsilon\eta)\cdot\widetilde{\bfv}=g(\sigma_s) \quad
  \mbox{in}\;\; \Bbb B^3,
\end{equation}
\begin{equation}
  -{\mathcal A}(\epsilon\eta)\widetilde{\bfv}
  +\vec{\mathcal B}(\epsilon\eta)\widetilde{p}
  ={1\over3}\vec{\mathcal B}(\epsilon\eta)(g(\sigma_s))
  \quad\mbox{in}\;\; \Bbb B^3,
\end{equation}
\begin{equation}
  \widetilde{\bfT}_{\epsilon\eta}(\widetilde{\bfv},\widetilde{p})
  \widetilde{\bfn}_{\epsilon\eta}=-\gamma\widetilde\kappa_{\epsilon\eta}
  \widetilde\bfn_{\epsilon\eta} \quad \mbox{on}\;\; \Bbb S^2,
\end{equation}
\begin{equation}
  \int_{|x|<1}\widetilde{\bfv}(x)G_{\epsilon\eta}(x)dx=0,
\end{equation}
\begin{equation}
  \int_{|x|<1}\widetilde{\bfv}(x)\cdot\widetilde{\bfw}_j^{\epsilon\eta}(x)
  G_{\epsilon\eta}(x)dx=0,\quad j=1,2,3.
\end{equation}
  We note that the above problem does have a unique solution. Indeed, from the
  proof of (2.34) of \cite{WuCui} we see that for any sufficiently small
  $\epsilon>0$, the conditions (2.24) and (2.25) are satisfied by $\varphi=
  g(\sigma_s)$, $\bfg={1\over3}\vec{\mathcal B}(\epsilon\eta)(g(\sigma_s))$ and
  $\bfh=-\gamma\widetilde\kappa_{\epsilon\eta}\widetilde\bfn_{\epsilon\eta}$ with
  $\rho=\epsilon\eta$. Hence the desired assertion follows from Lemma 2.2.

  Clearly, $\lim_{\epsilon\to 0}\widetilde{\bfv}=\bfv_s$, $\lim_{\epsilon\to 0}
  \widetilde{p}=p_s$ and $\bfV=\lim_{\epsilon\to 0}\epsilon^{-1}(\widetilde{\bfv}-
  \bfv_s)$. Hence, by a similar argument as in the proof of (3.6) and (3.7) we
  get, from (3.21) and (3.22) respectively, that
\begin{equation}
  \nabla\cdot\bfV=
  -[\vec{\mathcal B}'(0)\eta]\cdot\bfv_s \quad \mbox{in}\;\; \Bbb B^3,
\end{equation}
  and
\begin{equation}
  -\Delta\bfV +\nabla P
  =[\mathcal A'(0)\eta]\bfv_s-[\vec{\mathcal B}'(0)\eta]p_s+
  {1\over3}[\vec{\mathcal B}'(0)\eta]g(\sigma_s)
  \quad \mbox{in}\;\; \Bbb B^3,
\end{equation}
  where $P=\lim_{\epsilon\to 0}\epsilon^{-1}(\widetilde{p}-p_s)$. Next,
  recalling that
$$
  \widetilde{\bfn}_\rho(x)=\bfn(\phi_\rho(x)), \quad
  \widetilde{\kappa}_\rho(x)=\kappa(\phi_\rho(x)) \quad
  \mbox{for}\;\; x\in {\Bbb S}^2,
$$
  where as before $\bfn$ and $\kappa$ are the unit outward normal and the mean
  curvature of $\Gamma_\rho(=\phi_\rho(\Bbb S^2))$, respectively, by a direct
  computation we easily obtain
$$
  \widetilde\bfn_{\epsilon\eta}=\bfn_0-\epsilon\nabla_\omega\eta+o(\epsilon)
  \qquad \mbox{and}\qquad
  \widetilde\kappa_{\epsilon\eta}=1-\epsilon[\eta+{1\over2}
  \Delta_\omega\eta]+o(\epsilon).
$$
  Thus similarly as in the proof of (3.8) we get from (3.23) that
\begin{eqnarray*}
  [\nabla\otimes\bfV +\big(\nabla\otimes\bfV\big)^T\big]\bfn_0
  &=&-\Big\{[\vec{\mathcal B}'(0)\eta]\otimes\bfv_s+\big([\vec{\mathcal B}'(0)
  \eta]\otimes\bfv_s\big)^T\Big\}\bfn_0
  \nonumber\\
  &&+\gamma\{\nabla_\omega\eta+
  [\eta+{1\over2}\Delta_\omega\eta]\bfn_0\}+P\bfn_0
  +\bfT(\bfv_s,p_s)\nabla_\omega\eta \quad \mbox{on}\;\; \Bbb S^2.
\end{eqnarray*}
  A direct computation shows that (cf. (4.33) of \cite{FriedHu1})
$$
  \bfT(\bfv_s,p_s)\Big|_{\Bbb S^2}\nabla_\omega\eta
  =-(\gamma+2g(1))\nabla_\omega\eta.
$$
  Hence by using (1.11), (3.26) and the above result we obtain
\begin{eqnarray}
  \bfT(\bfV,P)\bfn_0
  &=&\displaystyle\;[\nabla\otimes\bfV +\big(\nabla\otimes\bfV\big)^T\big]
  \bfn_0-[P+\frac{2}{3}\nabla\cdot\bfV]\bfn_0
  \nonumber\\
  &=&\;-\Big\{[\vec{\mathcal B}'(0)\eta]\otimes\bfv_s
  +\big([\vec{\mathcal B}'(0)\eta]\otimes\bfv_s\big)^T\Big\}\bfn_0
  +{2\over3}\{[\vec{\mathcal B}'(0)\eta]\cdot\bfv_s\}\bfn_0
  \nonumber\\
  &&\;+\gamma[\eta+{1\over2}\Delta_\omega\eta]\bfn_0
  -2g(1)\nabla_\omega\eta \quad \mbox{on}\;\; \Bbb S^2.
\end{eqnarray}
  Finally, similarly as in the proof of (3.9) and (3.10) we get from (3.24)
  and (3.25) that
\begin{equation}
  \int_{|x|<1}[D_\rho\vec{\mathcal V}(\sigma_s,0)\eta
  +U(\eta)\bfv_s]\,dx=0,
\end{equation}
\begin{equation}
  \int_{|x|<1}\big\{[D_\rho\vec{\mathcal V}(\sigma_s,0)\eta]\cdot\bfw_j
  +\bfv_s\cdot\bfW_j(\eta)\big\}\,dx=0,\quad j=1,2,3,
\end{equation}
  respectively. Now let $\bfV_1=\bfV-v_s'(r)M(\eta)\bfn_0$ and $P_1=P-p_s'(r)M(\eta)$.
  Then from (3.6)--(3.11) and (3.26)--(3.30) we easily obtain
\begin{equation}
  \nabla\cdot\bfV_1=-g'(\sigma_s)\sigma_s'(r)M(\eta)\qquad \mbox{in}\;\; \Bbb B^3,
\end{equation}
\begin{equation}
  -\Delta\bfV_1+\nabla P_1=-{1\over3}\nabla[g'(\sigma_s)\sigma_s'(r)M(\eta)]
  \qquad\mbox{in}\;\; \Bbb B^3,
\end{equation}
\begin{equation}
  \bfT(\bfV_1,P_1)\bfn_0=\gamma(\eta+{1\over2}\Delta_\omega\eta)\bfn_0-
  2g(1)\nabla_\omega\eta+4g(1)\eta\bfn_0  \qquad\mbox{on}\;\;\Bbb S^2,
\end{equation}
\begin{equation}
  \int_{|x|<1}\bfV_1\,dx=0,
\end{equation}
\begin{equation}
  \int_{|x|<1}\bfV_1\cdot\bfw_j\,dx=0,\quad j=1,2,3.
\end{equation}
  In getting (3.35) we also used the fact that $\bfn_0\cdot\bfw_j=0$ for
  $j=1,2,3$. Using the relations $\displaystyle\int_{|x|=1}\!\!\nabla_\omega
  \eta\cdot\bfw_j\,dS_\omega=-\int_{|x|=1}\!\!\eta\nabla_\omega\cdot\bfw_j\,
  dS_\omega=0$ ($j=1,2,3$) we easily see that the relations (2.24) and (2.25)
  with $\rho=0$ are satisfied by
$$
  \varphi=-g'(\sigma_s)\sigma_s'(r)M(\eta),\qquad
  \bfg=-{1\over3}\nabla[g'(\sigma_s)\sigma_s'(r)M(\eta)],
$$
  and
$$
  \bfh=\gamma(\eta+{1\over2}\Delta_\omega\eta)\bfn_0-
  2g(1)\nabla_\omega\eta+4g(1)\eta\bfn_0.
$$
  Hence by Lemma 2.2 we see that the problem (3.31)--(3.35) has a unique
  solution $(\bfV_1,P_1)$ and, in particular, $\bfV_1$ is given by
\begin{eqnarray}
  \bfV_1
  \;&=&\;-\vec{\mathcal P}(0)\big[g'(\sigma_s)\sigma_s'(r)M(\eta)\big]
  -{1\over3}\mathbf Q(0)\big\{\nabla[g'(\sigma_s)\sigma_s'(r)M(\eta)]\big\}
  \nonumber\\
  &&+{\mathbf R}(0)\big[\gamma(\eta+{1\over2}\Delta_\omega\eta)\bfn_0
  -2g(1)\nabla_\omega\eta+4g(1)\eta\bfn_0\big],
\end{eqnarray}
  from which (3.20) immediately follows. This completes the proof of Lemma 3.2.
  $\qquad$ $\Box$
\medskip

  We are now ready to compute all the Fr\'{e}chet derivatives appearing in the
  right-hand side of (3.1). First, by (2.31), (3.3) and (3.19) we have
\begin{equation}
  D_{\tilde\sigma}\mathcal G(\sigma_s,0)v=\vec{\mathcal D}(0)
  [D_{\widetilde\sigma}\vec{\mathcal V}(\sigma_s,0)v]
  ={\rm tr}_{\Bbb S^2}\big\{\vec{\mathcal P}(0)[g'(\sigma_s)v]
  +{1\over 3}{\mathbf Q}(0)[\nabla(g'(\sigma_s)v)]\big\}
  \cdot\bfn_0.
\end{equation}
  Next, by (3.9) and the facts that $\vec{\mathcal V}(\sigma_s,0)=\bfv_s$,
  ${\rm tr}_{\Bbb S^2}(\bfv_s)=0$ we have
$$
  [\vec{\mathcal D}'(0)\eta]\vec{\mathcal V}(\sigma_s,0)=
  -{\rm tr}_{\Bbb S^2}[\vec{\mathcal V}(\sigma_s,0)]\cdot\nabla_\omega\eta=
  -{\rm tr}_{\Bbb S^2}(\bfv_s)\cdot\nabla_\omega\eta=0.
$$
  Thus by (2.31), (3.3) and (3.20) we have
\begin{eqnarray}
  D_\rho\mathcal G(\sigma_s,0)\eta\;&=&\;[\vec{\mathcal D}'(0)\eta]
  \vec{\mathcal V}(\sigma_s,0)+\vec{\mathcal D}(0)[D_\rho\vec{\mathcal V}
  (\sigma_s,0)\eta]
  \nonumber\\
  &=&\;{\rm tr}_{\Bbb S^2}[D_\rho\vec{\mathcal V}(\sigma_s,0)\eta]
  \cdot\bfn_0
  \nonumber\\
  &=&\;g(1)\eta-{\rm tr}_{\Bbb S^2}\Big\{\vec{\mathcal P}(0)
  \big[g'(\sigma_s)\sigma_s'(r)M(\eta)\big]+{1\over3}\mathbf Q(0)
  \big\{\nabla[g'(\sigma_s)\sigma_s'(r)M(\eta)]\big\}
  \nonumber\\
  &&\;-{\mathbf R}(0)\big[\gamma(\eta+{1\over2}
  \Delta_\omega\eta)\bfn_0-2g(1)\nabla_\omega\eta+4g(1)\eta
  \bfn_0\big]\Big\}\cdot\bfn_0.
\end{eqnarray}
  Thirdly, from (2.30) and a direct computation we have
\begin{eqnarray}
  D_{\tilde\sigma}\mathcal F(\sigma_s,0)v\;&=&\;\epsln^{-1}\mathcal A(0)v
  -\epsln^{-1}f'(\sigma_s)v+\chi(r-1)\Pi_1\big\{\vec{\mathcal D}(0)
  [D_{\widetilde\sigma}\vec{\mathcal V}(\sigma_s,0)v]\big\}
  \vec{\mathcal B}(0)\sigma_s\cdot\bfe_r
  \nonumber\\
  &&+\chi(r-1)\Pi_1\big(\vec{\mathcal D}(0)\vec{\mathcal V}(\sigma_s,0)\big)
  \vec{\mathcal B}(0)v\cdot\bfe_r
  \nonumber\\
  &=&\;\epsln^{-1}[\Delta-f'(\sigma_s)]v+\chi(r-1)\Pi_1\big\{\vec{\mathcal D}(0)
  [D_{\widetilde\sigma}\vec{\mathcal V}(\sigma_s,0)v]\big\}\sigma_s'(r)
  \nonumber\\
  &&+\chi(r-1)\Pi_1[{\rm tr}_{\Bbb S^2}(\bfv_s)\cdot\bfn_0]\nabla v
  \cdot\bfe_r
  \nonumber\\
  &=&\;\epsln^{-1}[\Delta-f'(\sigma_s)]v+\chi(r-1)\sigma_s'(r)
  \Pi_1\Big\{{\rm tr}_{\Bbb S^2}\big\{\vec{\mathcal P}(0)[g'(\sigma_s)v]
  \nonumber\\
  &&\;+{1\over 3}{\mathbf Q}(0)[\nabla(g'(\sigma_s)v)]\big\}
  \cdot\bfn_0\Big\}.
\end{eqnarray}
  Finally, from (2.30), (3.5), (3.38) and a direct computation we have
\begin{eqnarray}
  D_\rho\mathcal F(\sigma_s,0)\eta\;&=&\;\epsln^{-1}[\mathcal A'(0)\eta]\sigma_s
  +\chi(r-1)\Pi_1[D_\rho\mathcal G(\sigma_s,0)\eta]\vec{\mathcal B}(0)
  \sigma_s\cdot\bfe_r
  \nonumber\\
  &&
  +\chi(r-1)\Pi_1[\vec{\mathcal D}(0)\vec{\mathcal V}(\sigma_s,0)]
  [\vec{\mathcal B}'(0)\eta]\sigma_s\cdot\bfe_r
  \nonumber\\
  &=&\;\epsln^{-1}[\mathcal A'(0)\eta]\sigma_s+\chi(r-1)\Pi_1
  [D_\rho\mathcal G(\sigma_s,0)\eta]\sigma_s'(r)
  \nonumber\\
  && \qquad (\mbox{because}\;\;\vec{\mathcal D}(0)\vec{\mathcal V}(\sigma_s,0)
  ={\rm tr}_{\Bbb S^2}(\bfv_s)\cdot\bfn_0=0)
  \nonumber\\
  &=&-\epsln^{-1}[\Delta-f'(\sigma_s)]
  [\chi(r-1)\sigma_s'(r)\Pi_1(\eta)]+g(1)\sigma_s'(r)M(\eta)
  \nonumber\\
  &&-\chi(r-1)\sigma_s'(r)\Pi_1\Big\{{\rm tr}_{\Bbb S^2}\Big[\vec{\mathcal P}(0)
  \big[g'(\sigma_s)\sigma_s'(r)M(\eta)\big]
  \nonumber\\
  &&+{1\over3}\mathbf Q(0)\big\{\nabla[g'(\sigma_s)\sigma_s'(r)M(\eta)]\big\}
  -{\mathbf R}(0)\big\{\gamma(\eta+{1\over2}\Delta_\omega\eta)\bfn_0
  \nonumber\\
  &&\;-2g(1)\nabla_\omega\eta+4g(1)\eta
  \bfn_0\big\}\Big]\cdot\bfn_0\Big\}.
\end{eqnarray}
  In conclusion, we have the following lemma.
\medskip

  {\bf Lemma 3.3}\;\; {\em The Fr\'{e}chet derivative $D\Bbb F(U_s)$ is given
  by $(3.1)$, in which $D_{\widetilde\sigma}\mathcal F(\sigma_s,0)$,
  $D_\rho\mathcal F(\sigma_s,0)$, $D_{\widetilde\sigma}\mathcal G(\sigma_s,0)$
  and $D_\rho\mathcal G(\sigma_s,0)$ are given by $(3.39)$, $(3.40)$, $(3.37)$
  and $(3.38)$, respectively.} $\quad$ $\Box$
\medskip

  As usual, for a linear operator $L$ from a product space $X_1\times X_2$ to
  another product space $Y_1\times Y_2$ having the expression
$$
  L(x_1,x_2)=(L_{11}x_1+L_{12}x_2,L_{21}x_1+L_{22}x_2) \quad
  \mbox{for}\;\; (x_1,x_2)\in X_1\times X_2,
$$
  we write it as
$L=\left(
\begin{array}{cc}
L_{11}\; &\: L_{12}\\
L_{21}\; &\; L_{22}
\end{array}
\right)$.
  Then we have
$$
  D\Bbb F(U_s)=\left(
  \begin{array}{cc}
  D_{\widetilde\sigma}\mathcal F(\sigma_s,0) & \;\;D_\rho\mathcal F(\sigma_s,0)
  \\
  D_{\widetilde\sigma}\mathcal G(\sigma_s,0) & \;\; D_\rho\mathcal G(\sigma_s,0)
  \end{array}
  \right).
$$
  In the sequel we follow the idea of \cite{Cui2} to study the property of this
  operator and compute its spectrum.

  Recall that $m\ge3$, $m\in\Bbb N$ and $\theta\in
  (0,1)$. Let $\mathcal A_0: C^{m+\theta}(\overline{\Bbb B}^3)\to
  C^{m-2+\theta}(\overline{\Bbb B}^3)$ and $\mathcal J: C^{m+\theta}
  (\overline{\Bbb B}^3)\to C^{m+\theta}(\Bbb S^2)$ be the following operators:
\begin{equation}
  \mathcal A_0v=[\Delta-f'(\sigma_s)]v \quad
  \mbox{for}\;\;v\in C^{m+\theta}(\overline{\Bbb B}^3),
\end{equation}
\begin{equation}
  \mathcal J v={\rm tr}_{\Bbb S^2}\big[\vec{\mathcal P}(0)[g'(\sigma_s)v]+{1\over 3}
  {\mathbf Q}(0)(\nabla(g'(\sigma_s)v))\big]\cdot\bfn_0 \quad
  \mbox{for}\;\;v\in C^{m+\theta}(\overline{\Bbb B}^3).
\end{equation}
  Clearly, $\mathcal A_0\in L(C^{m+\theta}(\overline{\Bbb B}^3),
  C^{m-2+\theta}(\overline{\Bbb B}^3))$, and
\begin{equation}
  \mathcal J\in L(C^{m+\theta}(\overline{\Bbb B}^3),
  C^{m+\theta}(\Bbb S^2))\cap L(C^{m-2+\theta}(\overline{\Bbb B}^3),C^{m-2+\theta}
  (\Bbb S^2)).
\end{equation}
  Next let $\Pi_0: C^{m-2+\theta}(\Bbb S^2)\to C^{m-2+\theta}(\overline{\Bbb B}^3)$
  be the operator defined by $\Pi_0(\eta)=u$ for $\eta\in C^{m-2+\theta}(\Bbb S^2)$,
  where $u\in C^{m-2+\theta}(\overline{\Bbb B}^3)$ is the unique solution of the
  following boundary value problem:
$$
  [\Delta-f'(\sigma_s)]u=0\quad\mbox{in}\;\Bbb B^3,
  \quad u=\eta\quad\mbox{on}\;\Bbb S^2.
$$
  It is clear that $\Pi_0\in L(C^{m-2+\theta}(\Bbb S^2),
  C^{m-2+\theta}(\overline{\Bbb B}^3))\cap L(C^{m+\theta}(\Bbb S^2),
  C^{m+\theta}(\overline{\Bbb B}^3))$ and $\mathcal A_0\Pi_0=0$. We also let
  $\mathcal B_\gamma$ be the following operator from $C^{m+\theta}(\Bbb S^2)$ to
  $C^{m-1+\theta}(\Bbb S^2)$:
\begin{eqnarray}
  \mathcal B_\gamma\eta\;&=&\;g(1)\eta-\sigma_s'(1)\mathcal J\Pi_0(\eta)+
  {\rm tr}_{\Bbb S^2}\Big\{{\bf R}(0)\big[\gamma(\eta+{1\over2}
  \Delta_\omega\eta)\bfn_0
  \nonumber\\
  &&\;-2g(1)\nabla_\omega\eta+4g(1)\eta\bfn_0\big]\Big\}\cdot\bfn_0
  \quad\mbox{for}\;\;\eta\in C^{m+\theta}(\Bbb S^2).
\end{eqnarray}
  It is evident that $\mathcal B_\gamma\in L(C^{m+\theta}(\Bbb S^2),
  C^{m-1+\theta}(\Bbb S^2))$. Clearly,
\begin{equation}
\left\{
\begin{array}{rcl}
  D_{\widetilde\sigma}\mathcal F(\sigma_s,0)v\;&=&\;\epsln^{-1}\mathcal A_0v+
  \chi(r-1)\sigma_s'(r)\Pi_1{\mathcal J}v,\\
  D_{\widetilde\sigma}\mathcal G(\sigma_s,0)v\;&=&\;{\mathcal J}v,\\
  D_\rho\mathcal F(\sigma_s,0)\eta\;&=&\;\chi(r-1)\sigma_s'(r)
  \Pi_1(\eta)\big[\mathcal B_\gamma\eta+\sigma_s'(1)\mathcal J\Pi_0(\eta)
  -\mathcal J(\chi(r-1)\sigma_s'(r)\Pi_1(\eta))\big]\\
  &&\;-\epsln^{-1}\mathcal A_0[\chi(r-1)\sigma_s'(r)\Pi_1(\eta)],\\
  D_\rho\mathcal G(\sigma_s,0)\eta\;&=&\;\mathcal B_\gamma\eta+\sigma_s'(1)
  \mathcal J\Pi_0(\eta)-\mathcal J(\sigma_s'(r)\chi(r-1)\Pi_1(\eta)).
\end{array}
\right.
\end{equation}
  Finally, let $\Bbb M: \Bbb X_0\to \Bbb X$ and $\Bbb T: \Bbb X\to \Bbb X$ be
  the following operators:
$$
  \Bbb M=\left(
  \begin{array}{cc}
  \epsln^{-1}\mathcal A_0+\sigma_s'(1)\Pi_0\mathcal J
  \;\;&\;\; \sigma_s'(1)\Pi_0\mathcal B_\gamma
  \\
  \mathcal J\;\;&\;\; \mathcal B_\gamma
  \end{array}
  \right),
  \qquad
  \Bbb T=\left(
  \begin{array}{cc}
  I\;\;&\;\;-\sigma_s'(1)\Pi_0+\chi(r-1)\sigma_s'(r)\Pi_1
  \\
  0\;\;&\;\;I
  \end{array}
  \right).
$$
  It is easy to see that $\Bbb M\in L(\Bbb X_0,\Bbb X)$ and $\Bbb T\in L(\Bbb X)$.
  Moreover, since ${\rm tr}_{\Bbb S^2}\big\{-\sigma_s'(1)\Pi_0(\eta)+\chi(r-1)
  \sigma_s'(r)\Pi_1(\eta)\big\}=-\sigma_s'(1)\eta+\sigma_s'(1)\eta=0$ for any
  $\eta\in C^{m+\theta}(\Bbb S^2)$, we see that $\Bbb T$ maps $\Bbb X_0$ into
  itself, i.e. $\Bbb T\in L(\Bbb X)\cap L(\Bbb X_0)$. Besides, it can be easily
  seen that
$$
  \Bbb T^{-1}=\left(
  \begin{array}{cc}
  I\;\;&\;\;\sigma_s'(1)\Pi_0-\sigma_s'(r)\chi(r-1)\Pi_1
  \\
  0\;\;&\;\;I
  \end{array}
  \right).
$$
  By a simple computation we have
\medskip

  {\bf Lemma 3.4}\;\; $D\Bbb F(U_s)=\Bbb T\Bbb M\Bbb T^{-1}$. $\quad\Box$
\medskip

  Given a closed linear operator $L$ on a Banach space, we denote by
  $\sigma(L)$ and $\sigma_p(L)$ respectively the spectrum and the set of
  all eigenvalues of $L$. As an immediate consequence of Lemma 3.4 we have
  the following preliminary result:
\medskip

  {\bf Lemma 3.5}\;\;{\em $(i)$ Regarded as an unbounded linear operator in
  ${\Bbb X}$ with domain ${\Bbb X}_0$, $D\Bbb F(U_s)$ is an infinitesimal
  generator of an analytic semigroup in ${\Bbb X}$.

  $(ii)$ If $\delta$ is sufficiently small then for any $U\in \Bbb O_\delta$
  and in a neighborhood of $U_s$ in $\Bbb X_0$, we have that $D\Bbb F(U)$,
  regarded as an unbounded linear operator in ${\Bbb X}$ with domain ${\Bbb X}_0$,
  is an infinitesimal generator of an analytic semigroup in ${\Bbb X}$.

  $(iii)$ Let $V\in \Bbb X_0$ and $\lambda\in\Bbb C$. Then $D\Bbb F(U_s)V=
  \lambda V$ if and only if $\Bbb M(\Bbb T^{-1}V)=\lambda\Bbb T^{-1}V$.

  $(iv)$ $\sigma(D\Bbb F(U_s))=\sigma_p(D\Bbb F(U_s))=\sigma_p(\Bbb M)=
  \sigma(\Bbb M)$.}
\medskip

  {\em Proof:}\ \ $(i)$ By Lemma 3.4 it suffices to prove that the operator
  $\Bbb M$, regarded as an unbounded linear operator in ${\Bbb X}$ with
  domain ${\Bbb X}_0$, is an infinitesimal generator of an analytic semigroup
  in ${\Bbb X}$. We denote
$$
  \Bbb M_1=\left(
  \begin{array}{cc}
  \epsln^{-1}\mathcal A_0+\sigma_s'(1)\Pi_0\mathcal J
  \;\;&\;\; \sigma_s'(1)\Pi_0\mathcal B_\gamma
  \\
  0\;\;&\;\; \mathcal B_\gamma
  \end{array}
  \right),
  \qquad
  \Bbb M_2=\left(
  \begin{array}{cc}
  0 \;\;&\;\; 0
  \\
  \mathcal J\;\;&\;\;0
  \end{array}
  \right).
$$
  Then $\Bbb M=\Bbb M_1+\Bbb M_2$. From (3.43) it can be easily seen that
  $\Bbb M_2\in L(\Bbb X)$. Thus by a standard result for perturbations of
  generators of analytic semigroups (see \cite{Lunar}), we only need to show
  that $\Bbb M_1$, regarded as an unbounded linear operator in ${\Bbb X}$ with
  domain ${\Bbb X}_0$, is an infinitesimal generator of an analytic semigroup
  in ${\Bbb X}$.

  Clearly, the operator $\epsln^{-1}\mathcal A_0+\sigma_s'(1)\Pi_0\mathcal J$
  is an infinitesimal generator of an analytic semigroup in $C^{m-2+\theta}
  (\overline{\Bbb B}^3)$ (with domain $C^{m+\theta}(\overline{\Bbb B}^3)$).
  Next, from (3.42), (3.44) and the definition of $\Pi_0$ it can be easily seen
  that the operator $\mathcal B_\gamma$ can be rewritten in the following form:
$$
  \mathcal B_\gamma\eta={\rm tr}_{\Bbb S^2}(\vec{\upsilon})\cdot\bfn_0+g(1)\eta
  \quad\mbox{for}\;\;\eta\in C^{m+\theta}(\Bbb S^2),
$$
  where $\vec{\upsilon}$ is the second component of the solution $(\phi,
  \vec{\upsilon},\psi)$ of the following problem:
$$
  \Delta\phi=f'(\sigma_s)\phi
  \qquad\mbox{in}\;\;\Bbb B^3,
$$
$$
  \nabla\cdot\vec{\upsilon}=g'(\sigma_s)\phi
  \qquad\mbox{in}\;\;\Bbb B^3,
$$
$$
  -\Delta\vec{\upsilon}+\nabla\psi-
  {1\over3}\nabla(\nabla\cdot\vec{\upsilon}\,)=0
  \qquad\mbox{in}\;\;\Bbb B^3,
$$
$$
  \phi=-\sigma'_s(1)\eta
  \qquad\mbox{on}\;\;\Bbb S^2,
$$
$$
  \bfT(\vec{\upsilon},\psi){\bf n}_0=-2g(1)
  \nabla_\omega\eta+\gamma(\eta+{1\over2}
  \Delta_\omega\eta){\bf n}_0
  +4g(1)\eta\bfn_0\qquad\mbox{on}\;\;\Bbb S^2,
$$
$$
  \int_{|x|<1}\vec{\upsilon}\,dx=0,
$$
$$
  \int_{|x|<1}\vec{\upsilon}\times x\,dx=0.
$$
  This shows that $\mathcal B_\gamma$ is the same operator as that given by
  (3.9) of \cite{WuCui} with the same notation. Thus by Lemma 2.6 and Lemma
  3.1 of \cite{WuCui} we see that $\mathcal B_\gamma$ is an infinitesimal
  generator of an analytic semigroup in $C^{m-1+\theta}(\Bbb S^2)$ (with
  domain $C^{m+\theta}(\Bbb S^2)$). Besides, for any $\eta\in
  C^{m+\theta}(\Bbb S^2)$ we have
$$
  ||\sigma_s'(1)\Pi_0\mathcal B_\gamma\eta||_{C^{m-2+\theta}(\overline{\Bbb B}^3)}
  \le C||\mathcal B_\gamma\eta||_{C^{m-2+\theta}(\Bbb S^2)}
  \le C||\mathcal B_\gamma\eta||_{C^{m-1+\theta}(\Bbb S^2)},
$$
  i.e., $\sigma_s'(1)\Pi_0\mathcal B_\gamma$ is $\mathcal B_\gamma$-bounded in
  the notion of \cite{Nagel}. Hence, by Corollary 3.3 of \cite{Nagel} (see
  also Lemma 3.2 of \cite{Cui1}) we see that $\Bbb M_1$ is an infinitesimal
  generator of an analytic semigroup in ${\Bbb X}$ (with domain ${\Bbb X}_0$), as
  desired. This proves the assertion $(i)$.

  The assertion $(ii)$ is an easy consequence of the assertion $(i)$ and
  the fact that $\Bbb F\in C^\infty(\Bbb O_\delta,\Bbb X)$, and the assertion
  $(iii)$ is immediate. Finally, since $\Bbb X_0$ is compactly embedding into
  $\Bbb X$, the assertion $(iv)$ follows from $(i)$ and $(iii)$. The proof is
  complete. $\qquad$$\square$
\medskip

  Later on we always assume that $\delta$ is sufficiently small such that the
  open set $\Bbb O_\delta$ satisfies the condition of Lemma 3.4 $(ii)$.

  As in \cite{WuCui}, for every integer $l\geq 0$ we let $Y_{lm}(\omega)$,
  $m=-l,-l+1,\cdots,l-1,l$, be a normalized orthogonal basis (in $L^2(\Bbb S^2)$
  sense) of the space of all spherical harmonics of degree $l$. It is
  well-known that
$$
  \Delta_\omega Y_{lm}(\omega)=-(l^2+l)Y_{lm}(\omega).
$$
  For every such $l$ we denote by $F_{l}(r)$ the unique solution of the following
  ODE problem:
\begin{equation}
  \left\{
  \begin{array}{l}
  \displaystyle F_{l}''(r)+{2\over r}F_l'(r)-{l^2+l\over r^2}F_l(r)
  =f'(\sigma_s(r))F_{l}(r)
  \quad\mbox{for}\;\; 0<r<1,
  \\
  F_{l}'(0)=0,\qquad  F_{l}(1)=-\sigma'_s(1).
  \end{array}
  \right.
\end{equation}
  Next we denote
\begin{equation}
  \alpha_0=g(1)+\int_0^1g'(\sigma_s(r))F_0(r)r^2\,dr,
\end{equation}
\begin{equation}
  \alpha_l(\gamma)=-{l(l+2)(2l+1)\over 4(2l^2+4l+3)}(\gamma-\gamma_l),
  \quad \mbox{for}\;\; l\ge2,
\end{equation}
  where
\begin{equation}
  \gamma_l={4(2l+3)(l+1)\over l(l+2)(2l+1)}\Big[
  g(1)+\int_0^1g'(\sigma_s(r))F_l(r)r^{l+2}\,dr\Big], \quad  l\ge2.
\end{equation}
   From Lemma 3.2 and Lemma 3.3 of \cite{WuCui} we have:
\medskip

  {\bf Lemma 3.6}\hs {\em $(i)$ $\mathcal B_\gamma $ is a Fourier
  multiplication operator having the following expression: For any
  $\eta\in C^\infty({\Bbb S}^2)$ with Fourier expansion $\eta(\omega)=
  \sum_{l=0}^\infty\sum_{m=-l}^l c_{lm}Y_{lm}(\omega)$, we have
$$
  \mathcal B_\gamma\eta(\omega)=\alpha_0c_{00}Y_{00}+
  \sum_{l=2}^\infty\sum_{m=-l}^l\alpha_l(\gamma)c_{lm}
  Y_{lm}(\omega).
$$

  $(ii)$ The spectrum of $\mathcal B_\gamma$ is
  given by
$$
  \sigma(\mathcal B_\gamma)=\{\alpha_0,0\}\cup
  \{\alpha_l(\gamma):l=2,3,4,\cdots\}.
$$
  Moreover, the multiplicity of the eigenvalue $0$ is $3$.} $\quad$ $\Box$
\medskip

  By Lemma 3.4 $(ii)$ of \cite{WuCui} we know that $\gamma_l>0$ for all
  $l\geq 2$, and $\lim_{l\to\infty}\gamma_l=0$. Thus as in \cite{WuCui} we
  define
$$
  \gamma_*=\max_{l\geq 2}\gamma_l.
$$
  Clearly $0<\gamma_*<\infty$, and for $\gamma>\gamma_*$ we have
  $\alpha_l(\gamma)<0$ for all $l\geq 2$, while if $\gamma<\gamma_*$ then there
  exists $l\geq 2$ such that $\alpha_l(\gamma)>0$. Since clearly
  $\lim_{l\to\infty}\alpha_l(\gamma)=-\infty$, the following notation makes
  sense:
$$
  \alpha_\gamma^*=\max\{\alpha_0,\alpha_l(\gamma),l\ge2\}.
$$
  By Lemma 3.4 $(i)$ of \cite{WuCui} we know that $\alpha_0<0$. Thus
  $\alpha_\gamma^*<0$ for all $\gamma>\gamma_*$.

  In the following lemma $\epsln$ is the constant appearing in the equation
  (1.1), which also appears in the expressions of $D{\Bbb F}(U_s)$ and $\Bbb M$.
\medskip

  {\bf Lemma 3.7}\;\;{\em We have the following assertions:

  $(i)$ $0$ is an eigenvalue of $\Bbb M$ of multiplicity $3$.

  $(ii)$ For any $\gamma>0$ there exists a corresponding constant $\epsln_0'>0$
  and a bounded continuous function $\mu_{l,\gamma}$ defined on $(0,\epsln_0']$,
  such that for any $0<\epsln\le \epsln_0'$ we have
\begin{equation}
  \sigma(\Bbb M)\supseteq\{\lambda_{l,\gamma}(\epsln)\equiv\alpha_{l}(\gamma)+
  \epsln\mu_{l,\gamma}(\epsln),\;\;l=2,3,\cdots\},
\end{equation}
  and for each $l\geq 2$, the eigenvectors of $\Bbb M$ corresponding to the
  eigenvalue $\lambda_{l,\gamma}(\epsln)$ have the expression
$
  U_{lm}=\left(
  \begin{array}{c}
  \epsln a_{l,\gamma}(r,\epsln)\\
  1
  \end{array}
  \right)
  Y_{lm}(\omega),
$
  where $a_{l,\gamma}(r,\epsln)$ is a smooth function in
  $r\in [r_0,1]$ and is bounded continuous in $\epsln\in(0,\epsln_0']$.

  $(iii)$ For any $\gamma>\gamma_*$ there exists a corresponding constant
  $\epsln_0>0$ such that for any $0<\epsln\leq\epsln_0$, we have
\begin{equation}
  \sup\{{\rm Re}\lambda:\lambda\in\sigma(\Bbb M)\backslash\{0\}\}
  \leq {1\over2}\alpha_\gamma^*<0.
\end{equation}
}
\medskip

  {\em Proof}: We assert that for a vector $U=(v,\eta)\in\Bbb X_0$,
  ${\Bbb M}U=0$ if and only if $v=0$ and $\mathcal B_\gamma\eta=0$. Indeed, it
  is easy to see that ${\Bbb M}U=0$ if and only if ${\mathcal A}_0v=0$ and
  $\mathcal B_\gamma\eta=0$. Since $U\in\Bbb X_0$ implies that $v\in
  C^{m+\theta}(\overline{\Bbb B}^3)\cap C_0(\overline{\Bbb B}^3)$, we see that
  the boundary value of $v$ is zero. Hence, by the maximum principle we see
  that ${\mathcal A}_0v=0$ implies that $v=0$. This proves the desired
  assertion. By this assertion and the fact that $0$ is an eigenvalue of
  $\mathcal B_\gamma$ of multiplicity $3$, we immediately get the assertion
  $(i)$. Next, by making slight changes of the proof of Lemma 6.4 of
  \cite{Cui2}, we get the assertion $(ii)$. Finally, the assertion $(iii)$
  follows from a quite similar proof as that of Lemma 6.5 of \cite{Cui2}.
  $\qquad\square$
\medskip

\section{The proof of Theorem 1.1}
\setcounter{equation}{0}

  In this section we give the proof of Theorem 1.1. By Lemmas 2.1 and 2.3, we
  only need to prove that the stationary point $U_s$ of the equation (2.34) is
  asymptotically stable in case $\gamma>\gamma_*$ whereas unstable in case
  $\gamma<\gamma_*$. By Lemma 3.5 $(ii)$ we see that the equation (2.34) is
  of the parabolic type in a small neighborhood of $U_s$. It is thus natural to
  use the geometric theory for parabolic differential equations in Banach
  spaces to investigate this equation. In doing so, however, we meet with a
  serious difficulty that, by Lemma 3.5 $(iv)$ and Lemma 3.7 $(i)$, $0$ is an
  eigenvalue of $D{\Bbb F}(U_s)$, so that the standard linearized asymptotic
  stability principle cannot be used to tackle the equation (2.34). We shall as
  in \cite{Cui2} appeal to the Lie group action possessed by this equation to
  overcome this difficulty.

  For $\tau>0$ we denote by ${\mathbb B}^3_\tau$ the ball in ${\mathbb R}^3$
  centered at the origin with radius $\tau$. Regarding ${\mathbb B}^3_\tau$
  as a neighborhood of the unit element $0$ of the commutative Lie group
  ${\mathbb R}^3$, we see that $G={\mathbb B}^3_\tau$ is a local Lie
  group of dimension $3$. Given $z\in {\mathbb R}^3$, we denote by $S_z$
  the translation in ${\mathbb R}^3$ induced by $z$, i.e.,
$$
  S_z(x)=x+z \quad \mbox{for}\;\; x\in {\mathbb R}^3.
$$
  Let $\rho\in C^1({\mathbb S}^{2})$ such that $\|\rho\|_{C^1({\mathbb S}^{2})}$
  is sufficiently small, say, $\|\rho\|_{C^1({\mathbb S}^{2})}<\delta$ for some
  small $\delta>0$. For any $z\in{\mathbb B}^3_\tau$, consider the image of
  the hypersurface $r=1+\rho(\omega)$ under the translation $S_{z}$, which is
  still a hypersurface. This hypersurface has the equation
  $r=1+\tilde{\rho}(\omega)$ with $\tilde{\rho}\in C^1({\mathbb S}^{2})$, and
  $\tilde{\rho}$ is uniquely determined by $\rho$ and $z$. We denote
$$
  \tilde{\rho}=S_z^*(\rho).
$$
  By some similar arguments as in the proof of Lemmas 4.1 and 4.3 of \cite{Cui2}
  we can show that if $\rho\in O^{m+\theta}_\delta({\mathbb S}^{2})$
  then $S_z^*(\rho)\in C^{m+\theta}({\mathbb S}^{2})$ and $S_z^*\in
  C(O^{m+\theta}_\delta({\mathbb S}^{2}),C^{m+\theta}({\mathbb S}^{2}))
  \cap C^1(O^{m+\theta}_\delta({\mathbb S}^{2}),
  C^{m-1+\theta}({\mathbb S}^{2}))$. Next, for each $z\in {\mathbb B}^3_\tau$
  and $\rho\in C^1({\mathbb S}^{2})$ such that $\|\rho\|_{C^1({\mathbb S}^{2})}$
  is sufficiently small, let $P_{z,\rho}:C(\overline{\mathbb B}^3)\to
  C(\overline{\mathbb B}^3)$ be the mapping
$$
   P_{z,\rho}(u)(x)=u(\Phi_{\rho}^{-1}(\Phi_{S_z^*(\rho)}(x)-z))\quad
   \mbox{for}\;\;u\in C(\overline{\mathbb B}^3).
$$
  This mapping is well-defined and, actually, we have $P_{z,\rho}\in
  L(C(\overline{\mathbb B}^3))$. Indeed, letting $\tilde{\rho}=S_z^*(\rho)$
  and denoting by $\Omega_\rho$ and $\Omega_{\tilde{\rho}}$ the domains enclosed
  by the hypersurfaces $r=1+\rho(\omega)$ and $r=1+\tilde{\rho}(\omega)$,
  respectively, we see that $x\in\Bbb B^3$ if and only if
  $\Phi_{S_z^*(\rho)}(x)=\Phi_{\tilde{\rho}}(x)\in\Omega_{\tilde{\rho}}$,
  which is equivalent to $\Phi_{S_z^*(\rho)}(x)-z\in\Omega_\rho$. But we
  know that $y\in\Omega_\rho$ if and only if $\Phi_{\rho}^{-1}(y)\in\Bbb B^3$.
  Hence, we see that $x\in\Bbb B^3$ if and only if
  $\Phi_{\rho}^{-1}(\Phi_{S_z^*(\rho)}(x)-z)\in\Bbb B^3$. Moreover, it can
  be easily seen that for $\rho\in C^1({\mathbb S}^{2})$ the mapping $x\to
  \Phi_{\rho}^{-1}(\Phi_{S_z^*(\rho)}(x)-z)$ is a $C^1$-diffeomorphism from
  $\overline{\Bbb B}^3$ onto itself. Hence the desired assertion follows.
  Furthermore, it is also immediate to see that if $\rho\in O^{m+\theta}_\delta
  ({\mathbb S}^{2})$ then the mapping $x\to\Phi_{\rho}^{-1}(\Phi_{S_z^*(\rho)}(x)-z)$
  is a $C^{m+\theta}$-diffeomorphism from $\overline{\Bbb B}^3$ onto itself,
  so that $P_{z,\rho}\in L(C^{m+\theta}(\overline{\mathbb B}^3))$.

  Now, for each $z\in G={\mathbb B}^3_\tau$ we define a mapping
  ${\mathbf S}_z^*:\Bbb O_\delta\to\Bbb X_0$ (recall that $\Bbb O_\delta=
  (C^{m+\theta}(\overline{\Bbb B}^3)\cap C_0(\overline{\Bbb B}^3))
  \times O^{m+\theta}_\delta({\Bbb S}^2)$ and $\Bbb X_0=
  (C^{m+\theta}(\overline{\Bbb B}^3)\cap C_0(\overline{\Bbb B}^3))
  \times C^{m+\theta}({\Bbb S}^2)$) as follows: For any $U=(u,\rho)\in
  \Bbb O_\delta$, let
$$
  {\mathbf S}_z^*(u,\rho)=(P_{z,\rho}(u),S_z^*(\rho)).
$$
  It is obvious that if $u\in C_0(\overline{\Bbb B}^3))$ then also
  $P_{z,\rho}(u)\in C_0(\overline{\Bbb B}^3))$. Thus for $\tau$ and
  $\delta$ sufficiently small this mapping is well-defined and it
  does map $\Bbb O_\delta$ into $\Bbb X_0$. We have
\medskip

  {\bf Lemma 4.1}\ \ {\em Let $m\geq 3$ and $0<\theta<1$. Let
$$
  \Bbb O_\delta'=(C^{m-2+\theta}(\overline{\Bbb B}^3)\cap
  C_0(\overline{\Bbb B}^3))\times O^{m-1+\theta}_\delta({\Bbb S}^2)
  \subseteq\Bbb X.
$$
  For sufficiently small $\tau>0$ and $\delta>0$ we have the following
  assertions:

  $(i)$ For any $z\in {\mathbb B}^3_\tau$ we have ${\mathbf S}_z^*
  \in C(\Bbb O_\delta',\Bbb X)\cap C(\Bbb O_\delta,\Bbb X_0)$.

  $(ii)$ For any $z,w\in {\mathbb B}^3_\tau$ we have
$$
  {\mathbf S}_z^*\circ {\mathbf S}_w^*={\mathbf S}_{z+w}^*,\quad
  {\mathbf S}_0^*=id, \quad \mbox{and}
  \quad ({\mathbf S}_z^*)^{-1}={\mathbf S}_{-z}^*.
$$

  $(iii)$ The mapping ${\mathbf S}^*: z\to {\mathbf S}_z^*$ from
  ${\mathbb B}^3_\tau$ to $C(\Bbb O_\delta',\Bbb X)$ is an injection, and
\begin{equation}
  {\mathbf S}^*\in C^k({\mathbb B}^3_\tau,C^l(\Bbb O_\delta,
  C^{m-k-l+\theta}(\overline{\mathbb B}^3)\times C^{m-k-l+\theta}
  ({\mathbb S}^{2}))),\quad k\geq 0, \;\; l\geq 0,\;\; k+l\leq m.
\end{equation}
  Moreover, for fixed $z\in {\mathbb B}^3_\tau$ we have
\begin{equation}
  D{\mathbf S}^*_z\in C(\Bbb O_\delta,
  L(C^{m-1+\theta}(\overline{\mathbb B}^3)\times C^{m-1+\theta}
  ({\mathbb S}^{2}))),
\end{equation}
  i.e., for any $U\in\Bbb O_\delta$, the operator $D{\mathbf S}^*_z(U)$
  $($which is, by $(4.1)$, an bounded linear operator from
   $C^{m+\theta}(\overline{\mathbb B}^3)\times C^{m+\theta}
  ({\mathbb S}^{2})$ to $C^{m-1+\theta}(\overline{\mathbb B}^3)\times
  C^{m-1+\theta}({\mathbb S}^{2}))$ can be extended into a bounded
  linear operator from $C^{m-1+\theta}(\overline{\mathbb B}^3)\times
  C^{m-1+\theta}({\mathbb S}^{2})$ to itself, and the mapping $U\to
  D{\mathbf S}^*_z(U)$ from $\Bbb O_\delta$ to $L(C^{m-1+\theta}
  (\overline{\mathbb B}^3)\times C^{m-1+\theta}({\mathbb S}^{2}))$ is
  continuous.

  (iv) Define $p:{\mathbb B}^3_\tau\times \Bbb O_\delta'\to\Bbb X$ by
  $p(z,U)={\mathbf S}_z^*(U)$ for $(z,U)\in {\mathbb B}^3_\tau\times
  \Bbb O_\delta'$. Then for any $U\in \Bbb O_\delta$ we have
  $p(\cdot,U)\in C^1({\mathbb B}^3_\tau,\Bbb X)$, and ${\rm rank}\,D_z
  p(z,U)=3$ for every $z\in {\mathbb B}^3_\tau$ and $U\in \Bbb O_\delta$.
  If furthermore $U\in \Bbb Z=C^\infty(\overline{\Bbb B}^3)\times
  C^\infty({\mathbb S}^{2})$ then $p(\cdot,U)\in C^\infty({\mathbb
  B}^3_\epsln,\Bbb Z)$.}
\medskip

  {\em Proof}:\ \ This lemma follows from a similar argument as that in the
  establishment of Lemma 4.4 of \cite{Cui2}. Indeed, to get the assertions of
  this lemma we only need to replace the Sobolev and corresponding Besov
  spaces $W^{m,q}({\Bbb B}^3)$, $B^{m-1-1/q}({\Bbb S}^2)$ in \cite{Cui2}
  with the H\"{o}lder spaces $C^{m+\theta}(\overline{\Bbb B}^3)$ and
  $C^{m+\theta}({\Bbb S}^2)$, and after a such replacement all the analysis
  presented in Section 4 of \cite{Cui2} still holds. In particular, the proof
  of (4.1) uses a similar argument as that of the second assertion in $(i)$
  of Lemma 4.3 of \cite{Cui2}. Since this analysis is lengthy but does not have
  new ingredient different from that in \cite{Cui2}, we omit it here.
  $\quad\Box$
\medskip

  {\bf Corollary 4.2}\ \ {\em Let assumptions be as in {\rm Lemma} $4.1$.
  Then for sufficiently small $\tau>0$, $\delta>0$ and fixed $z\in
  {\mathbb B}^3_\tau$ we also have the following relation:}
\begin{equation}
  D{\mathbf S}^*_z\in C(\Bbb O_\delta,
  L(C^{m-2+\theta}(\overline{\mathbb B}^3)\times C^{m-1+\theta}
  ({\mathbb S}^{2}))).
\end{equation}

  {\em Proof}:\ \ Since ${\mathbf S}_z^*(u,\rho)=(P_{z,\rho}(u),S_z^*(\rho))$,
  we see that
$$
  D{\mathbf S}^*_z(u,\rho)(v,\eta)=(D_u P_{z,\rho}(u)v+D_\rho
  P_{z,\rho}(u)\eta,DS_z^*(\rho)\eta).
$$
  From (4.2) we have $[(u,\rho)\to DS_z^*(\rho)]\in C(\Bbb O_\delta,
  L(C^{m-1+\theta}({\mathbb S}^{2})))$ and $[(u,\rho)\to D_\rho
  P_{z,\rho}(u)]\in C(\Bbb O_\delta,L(C^{m-1+\theta}({\mathbb S}^{2})))
  \subseteq C(\Bbb O_\delta,L(C^{m-1+\theta}({\mathbb S}^{2}),
  C^{m-2+\theta}({\mathbb S}^{2})))$. Since $m-1\geq m-2$, by (4.2) we
  also have $[(u,\rho)\to D_u P_{z,\rho}(u)]\in C(\Bbb O_\delta,
  L(C^{m-2+\theta}(\overline{\mathbb B}^3))$. Combining these
  assertions together, we see that (4.3) follows. $\qquad\Box$
\medskip

  In the sequel, for $\rho=\rho(t)$, $u=u(x,t)$ and $U=(u(x,t),\rho(t))$, we
  denote by $P_{z,\rho}(u)$ the function $\widetilde{u}(x,t)=
  u(\Phi_{\rho(t)}^{-1}(\Phi_{S_z^\ast(\rho(t))}(x)-z),t)$, by
  $S_z^\ast(\rho)$ the function $\widetilde{\rho}(t)=S_z^\ast(\rho(t))$, and
  by ${\mathbf S}_z^\ast (U)$ the vector function $(P_{z,\rho}(u),
  S_z^\ast(\rho))=(\widetilde{u}(x,t),\widetilde{\rho}(t))$.
\medskip

  {\bf Lemma 4.3}\ \ {\em If $U=(u,\rho)$ is a solution of the equation
  ${dU}\!/{dt}={\Bbb F}(U)$ such that $\|\rho\|_{C^1({\mathbb S}^2)}$ is
  sufficiently small, then for any $z\in {\mathbb R}^3$ such that $|z|$ is
  sufficiently small, ${\mathbf S}_z^*(U)=(P_{z,\rho}(u),S_z^\ast(\rho))$ is
  also a solution this equation.}
\medskip

  {\em Proof}:\ \ It is easy to see that if $(\sigma,\bfv,p,\Omega)$ is a
  solution of the system of equations (1.1)--(1.8), then for any $z\in
  {\mathbb R}^3$, we have that $(\widetilde{\sigma},\widetilde{\bfv},
  \widetilde{p},\widetilde{\Omega})$ defined by
$$
  \widetilde{\sigma}(x,t)=\sigma(x-z,t), \quad \widetilde{p}(x,t)=p(x-z,t),
  \quad \widetilde{\bfv}(x,t)=\bfv(x-z,t), \quad
  \widetilde{\Omega}(t)=\Omega(t)+z,
$$
  is also a solution of (1.1)--(1.8). From this fact we see immediately that
  if $U=(u,\rho)$ is a solution of the equation
\begin{equation}
  \frac{dU}{dt}={\Bbb F}(U),
\end{equation}
  then $\widetilde{U}=(\widetilde{u},\widetilde{\rho})$, where
$$
  \widetilde{u}(x,t)=u(\Phi_{\rho(t)}^{-1}
  (\Phi_{S_z^\ast(\rho(t))}(x)-z),t),
  \quad \widetilde{\rho}(t)=S_z^\ast(\rho(t)),
$$
  is also a solution of this equation, which is the desired assertion.
  $\qquad\Box$
\medskip

  {\bf Lemma 4.4}\ \ {\em If $\tau$ and $\delta$ are sufficiently small then
  for any $z\in {\mathbb B}^3_\tau$ and $U=(u,\rho)\in {\Bbb O}_\delta$ we
  have}
\begin{equation}
  {\mathbb F}({\mathbf S}_z^\ast (U))=D{\mathbf S}_{z}^*(U){\mathbb F}(U).
\end{equation}

  {\em Proof}:\ \ By Lemma 3.5 $(ii)$ we see that the equation (4.1) is of the
  parabolic type in ${\Bbb O}_\delta$, provided $\delta$ is sufficiently small.
  Hence, by a well-known result in the theory of differential equations of the
  parabolic type in Banach spaces (cf. Theorem 8.1.1 of \cite{Lunar}) we see
  that given any $U_0=(u,\rho)\in\Bbb X_0$ there exists $t_0>0$ such that the
  equation (4.1) has a unique solution $U=U(t)$ for $0\leq t\leq t_0$, which
  belongs to $C([0,t_0],\Bbb X)\cap C((0,t_0],\Bbb O_\delta)\cap L^\infty
  ((0,t_0),\Bbb X_0)\cap C^1((0,t_0],\Bbb X)$ and satisfies the initial
  condition $U(0)=U_0$. Let $\widetilde{U}(t)={\mathbf S}_z^*(U(t))$ for $0\leq
  t\leq t_0$. By Lemma 4.2, $\widetilde{U}$ is also a solution of (4.1),
  satisfying the initial condition $\widetilde{U}(0)={\mathbf S}_z^\ast(U_0)$.
  The fact that $\widetilde{U}$ is the solution of (4.1) implies that
$$
  {d\widetilde{U}(t)\over dt}={\mathbb F}(\widetilde{U}(t)) \quad
  \mbox{for}\;\; 0<t\leq t_0.
$$
  On the other hand, since $\widetilde{U}(t)={\mathbf S}_z^\ast (U(t))$, we
  have
$$
  {d\widetilde{U}(t)\over dt}=D{\mathbf S}_z^\ast (U(t)){dU(t)\over dt}
  =D{\mathbf S}_z^\ast (U(t)){\mathbb F}(U(t))
   \quad \mbox{for}\;\; 0<t\leq t_0.
$$
  Thus
\begin{equation}
  {\mathbb F}(\widetilde{U}(t))=D{\mathbf S}_z^\ast (U(t)){\mathbb F}
  (U(t))\quad \mbox{for}\;\; 0<t\leq t_0.
\end{equation}
  If $U(t)$ is a strict solution (in the sense of \cite{Lunar}) then $U\in
  C([0,t_0],\Bbb O_\delta)\cap C^1([0,t_0],\Bbb X)$ and clearly
  $\widetilde{U}(t)$ is also a strict solution, so that by directly letting
  $t\to0^+$ we get
\begin{equation}
  {\mathbb F}({\mathbf S}_z^\ast (U_0))=
  D{\mathbf S}_{z}^*(U_0){\mathbb F}(U_0).
\end{equation}
  If $U(t)$ is not a strict solution then we establish this relation in the
  following way: Let
$$
   \widetilde{\Bbb X}_0=C^{m-2+\theta}(\overline{\Bbb B}^3)\times C^{m-2+\theta}
  ({\Bbb S}^2), \quad
  \widetilde{\Bbb X}=C^{m-4+\theta}(\overline{\Bbb B}^3)\times C^{m-3+\theta}
  ({\Bbb S}^2),
$$
  and let $\widetilde{\Bbb O}_\delta$ be a small neighborhood of the
  origin of $\widetilde{\Bbb X}_0$. Then $\Bbb F\in C(\widetilde{\Bbb O}_\delta,
  \widetilde{\Bbb X})$ and $D{\mathbf S}_{z}^*\in C(\widetilde{\Bbb O}_\delta,
  L(\widetilde{\Bbb X}))$ (by (4.3)). Hence, since $U\in C([0,t_0],\Bbb
  X)\subseteq C([0,t_0],\widetilde{\Bbb X}_0)$, which also ensures that
  $\widetilde{U}\in C([0,t_0],\widetilde{\Bbb X}_0)$,
  by letting $t\to0^+$ we see that (4.7) still holds.   $\qquad\Box$
\medskip

  We are now ready to prove Theorem 1.1.
\medskip

  {\bf Proof of Theorem 1.1}:\ \ We first assume that $\gamma>\gamma_*$. Then
  by Lemma 3.5 $(iv)$ and Lemma 3.7 $(i)$ and $(iii)$ we see that
$$
  \sup\{{\rm Re}\lambda:\lambda\in\sigma(D{\Bbb F}(U_s))\backslash\{0\}\}
  \leq {1\over2}\alpha_\gamma^*<0,
$$
  and $0$ is an eigenvalue of $D{\Bbb F}(U_s)$ of multiplicity $3$. By Lemma
  4.1 we see that the mapping $p:{\mathbb B}^3_\tau\times \Bbb O_\delta'\to
  \Bbb X$ defined by $p(z,U)={\mathbf S}_z^*(U)$ for $(z,U)\in {\mathbb B}^3_\tau
  \times\Bbb O_\delta'$ is a Lie group action of the local Lie group $G=
  {\mathbb B}^3_\tau$ to the space $\Bbb X$, and by Lemma 4.4 we see that the
  equation $dU\!/dt={\Bbb F}(U)$ is quasi-invariant (in the sense of
  \cite{Cui2}) under this Lie group action. Since for every $z\in
  {\mathbb B}^3_\tau$ and $U\in \Bbb O_\delta$ we have ${\rm rank}\,
  D_z p(z,U)=3=$ the multiplicity of the eigenvalue $0$ of $D{\Bbb F}(U_s)$, we
  see that all the conditions of Theorem 2.1 of \cite{Cui2} are satisfied by
  the equation $dU\!/dt={\Bbb F}(U)$ and the Lie group action $(G,p)$. Hence,
  by Theorem 2.1 of \cite{Cui2} we conclude that there exist two neighborhoods
  ${\mathbb O}_1$ and ${\mathbb O}_2$ of $U_s$ in $\Bbb X_0$, ${\mathbb O}_1$
  is larger than ${\mathbb O}_2$ (i.e. ${\mathbb O}_2\subseteq {\mathbb O}_1$),
  such that for each $U_0\in {\mathbb O}_1$ the initial value problem (2.34)
  has a unique solution $U=U(t)$ for all $t\geq 0$, which belongs to
  $C([0,\infty),\Bbb X)\cap C((0,\infty),{\mathbb O}_1)\cap L^\infty
  ((0,\infty),\Bbb X_0)\cap C^1((0,\infty),\Bbb X)$, and for any $U_0\in
  {\mathbb O}_2$ there exists corresponding $V_0\in {\mathbb O}_1$ and $z\in
  G={\mathbb B}^3_\tau$ which are uniquely determined by $U_0$, such that
$$
  U_0={\mathbf S}_z^*(V_0), \quad \lim_{t\to\infty}U(t;V_0)=U_s, \quad
  \mbox{and} \quad \lim_{t\to\infty}U(t;U_0)={\mathbf S}_z^*(U_s),
$$
  where $U(t;V_0)$ refers to the solution of (2.34) with initial value $V_0$,
  and similarly for $U(t;U_0)$. Moreover, the convergence in the above limit
  relations is exponentially fast. From this result and Lemmas 2.1 and 2.3 we
  immediately obtain (1.13).

  Next we assume that $\gamma<\gamma_*$. Then there exists $l\geq 2$ such that
  $\alpha_l(\gamma)>0$. By Lemma 3.5 $(iv)$ and Lemma 3.7 $(ii)$ we see that
  for $\epsln>0$ sufficiently small, $D{\Bbb F}(U_s)$ has an eigenvalue
  $\lambda_{l,\gamma}(\epsln)=\alpha_l(\gamma)+\epsln\mu_{l,\gamma}(\epsln)$
  which is clearly positive. Hence, by a standard result in the geometric
  theory of parabolic differential equations in Banach spaces (cf. Theorem
  9.1.3 of \cite{Lunar}) we conclude that $U_s$ is unstable as a stationary
  point of the equation $dU\!/dt={\Bbb F}(U)$. Using again Lemmas 2.1 and 2.3,
  we get the last assertion of Theorem 1.1. This completes the proof.
  $\qquad\Box$
\medskip

   {\bf Acknowledgement}.\hskip 1em This work is supported by the National
   Natural Science Foundation of China under the grant number 10771223 and
   a funds in Sun Yat-Sen University.

\vsss

{\small
\begin{thebibliography}{99}

\bibitem{Cui1} S. Cui, Well-posedness of a multidimensional free
  boundary problem modelling the growth of nonecrotic tumors, {\em J.
  Func. Anal.}, 245(2007), 1--18.

\bibitem{Cui2} S. Cui, Lie group action and stability analysis of
  stationary solutions for a free boundary problem modelling tumor growth,
  preprint (see arXiv: 0712.2483vl).

\bibitem{EscSim} J. Escher and G. Simonett,  A center monifold analysis for the
  Mullins-Sekerka model, {\em J. Diff. Equa.}, 143(1998), 267--292.

\bibitem{Franks1} S. J. H. Franks, H. M. Byrne, J. P. King, J. C. E. Underwood,
  C. E. Lewis,  Modelling the early growth of ductal carcinoma in situ of
  the breast, {\em J. Math. Biol.}, 47(2003), 424--452.

\bibitem{Franks2} S. J. H. Franks, H. M. Byrne, J. P. King, J. C. E. Underwood,
  C. E. Lewis,  Modelling the growth of comedo ductal carcinoma in situ,
  {\em Math. Med. Biol.}, 20(2003), 277--308.

\bibitem{Franks3} S. J. H. Franks, H. M. Byrne, J. C. E. Underwood,
  C. E. Lewis,  Biological inferences from a mathematical model of comedo
  ductal carcinoma in situ of the breast, {\em J. Theoret. Biol.},
  232(2005), 523--543.

\bibitem{Franks4} S. J. H. Franks, J. P. King,  Interactions between
  a uniformly proliferating tumour and its surroundings: Uniform material
  properties, {\em Math. Med. Biol.}, 20(2003), 47--89.

\bibitem{Fried1} A. Friedman,  A free boundary problem for a coupled system
  of elliptic, hyperbolic, and Stokes equations modeling tumor growth,
  {\em Interfaces and Free Boundaries}, 8(2006), 247--261.

\bibitem{Fried2} A. Friedman,  Mathematical analysis and challenges arising
  from models of tumor growth, {\em Mathematical Models and Methods in
  Applied Sciences}, 17, suppl.(2007), 1751--1772.

\bibitem{FriedHu1} A. Friedman and B. Hu, Bifurcation for a free boundary
  problem modeling tumor growth by Stokes equation, {\em SIAM J. Math. Anal.},
  39(2007), 174--194.

\bibitem{FriedHu2} A. Friedman and B. Hu, Bifurcation from stability to
  instability for a free boundary problem modeling tumor growth by Stokes
  equation, {\em J. Math. Anal. Appl.}, 327(2007), 643--664.

\bibitem{Lunar} A. Lunardi,  Analytic Semigroups and Optimal Regularity in
  Parabolic Problems,  Birkh$\ddot{\text{a}}$user,  Basel, 1995.

\bibitem{Nagel} R. Nagel,  Towards a ``Matrix Theory'' for Unbounded Operator
Matrices,  {\em Math. Z.}, 201(1989), 57--68.

\bibitem{WuCui} J. Wu and S. Cui,  Asymptotic behavior of solutions of a
  free boundary problem modelling the growth of tumors with Stokes equations,
  preprint (see arXiv: 0806.1353vl).

\end{thebibliography}
}
\end{document}